\documentclass[12pt]{article}
\usepackage[dvips]{graphicx}
\usepackage{amsmath,amsfonts,amssymb,amsthm,color}

\usepackage{bm}

\newtheorem{theorem}{Theorem}
\newtheorem{lemma}{Lemma}
\newtheorem{corollary}{Corollary}
\newtheorem{remark}{Remark}

\newcommand{\bdeltap}{\bm{\delta}_{\phi}}
\newcommand{\bdelta}{\bm{\delta}}

\title{Improved estimation of the MSEs and \\
the MSE matrices for shrinkage estimators of multivariate normal means
and their applications} 
\author{Hisayuki Hara\\
Department of Geosystem Engineering\\
University of Tokyo}
\date{October 2007}
\begin{document}
\maketitle
\begin{abstract}
 In this article we provide some nonnegative and positive estimators 
 of the mean squared errors(MSEs) for shrinkage estimators of
 multivariate normal means. Proposed estimators are shown to improve on
 the uniformly minimum variance unbiased estimator(UMVUE) under
 a quadratic loss criterion.      
 A similar improvement is also obtained for the estimators of the MSE
 matrices for shrinkage estimators. 
 We also apply the proposed estimators of the MSE matrix to form
 confidence sets centered at shrinkage estimators and show their
 usefulness through numerical experiments.
\end{abstract}

\section{Introduction.}
\label{sec:1}
\hspace*{12pt}
Since Stein \cite{Stein-1956} and James and Stein \cite{James-Stein}
proved the inadmissibility of the sample mean as an estimator of a 
multivariate normal mean and demonstrated the superiority of the
James-Stein estimator, considerable efforts have been devoted to
develop the theory of the shrinkage estimation in theoretical
literatures.  
In step with the development, shrinkage estimators have been applied 
to some practical problems. 
For example, see Fay and Herriot\cite{Fay-Herriot}, Battese, Harter
and Fuller\cite{Battese-etal} for the small-area estimation and 
Adkins and Eells\cite{Adkins-Eells} for the estimation of energy
demand. 
As pointed out in Wan et al.\cite{Wan-etal}, however, shrinkage
estimation does not necessarily find many applications in empirical
work. This is considered mainly because relatively little attention has
been paid to evaluate the precision of shrinkage estimators.  
This issue is important not only for recognizing how much it pays 
to use shrinkage estimators for a given problem but also for
constructing confidence sets. 

Joshi\cite{Joshi} proved the existence of confidence sets of a
multivariate normal mean centered at shrinkage estimators
which show higher coverage probabilities than the conventional one, a
sphere centered at the observations, when the variance is known.  
Hwang and Casella\cite{Hwang-Casella} give a explicit proof that when
the variance is known, a confidence set centered at the positive-part
Stein estimator has uniformly higher coverage probability than the
conventional one. 

Carter et al.\cite{Carter-etal} discuss this problem with unknown
variance in the context of linear regression models and proposed to
use the estimators of the MSE and MSE matrix of the James-Stein
estimator of regression coefficients to evaluate the precision of it. 
They derived the UMVUEs of the MSE and the MSE matrix for the James-Stein 
estimator and proposed to use the UMVUE of the MSE matrix to form a
confidence set. 
They also showed that the proposed confidence set has an asymptotically
smaller volume than that of the conventional one based on $F$
statistics.  
Wan et al.\cite{Wan-etal} generalized their argument to a wider class of  
shrinkage estimators.  
As clarified in Section 2, however, the UMVUEs of the MSE and the MSE
matrix of a shrinkage estimator have a serious drawback that they
cannot always take positive or positive definite.  Hence in a practical
sense, the results in Carter et al.\cite{Carter-etal} and Wan et
al.\cite{Wan-etal} may still not be suitable for empirical work. 

Kubokawa and Srivastava\cite{Kubokawa-Srivastava} focused their
attention on the drawback of the UMVUEs and considered the estimation of 
the MSEs and the MSE matrices from a decision theoretical viewpoint.
They defined quadratic loss functions for evaluating estimators and
proposed positive estimators of the MSE and the MSE matrix of the 
James-Stein estimator improving on the UMVUE under the loss functions.

In this article we take a similar approach to the one in Kubokawa and
Srivastava\cite{Kubokawa-Srivastava} and extend their argument to the
MSEs and the MSE matrices of a wider class of shrinkage estimators.   
We propose nonnegative and positive estimators of the MSEs and the 
MSE matrices of shrinkage estimators improving on the UMVUEs. 
We also apply the proposed estimators of the MSE matrices to form 
confidence sets centered at shrinkage estimators and show that they 
have higher coverage probabilities than the conventional one through
Monte Carlo studies. 

The organization of this article is as follows. 
In Section 2 we summarize some preliminary facts on shrinkage 
estimators and their MSEs and MSE matrices.
Section 3 and Section 4 provide some nonnegative and positive estimators
of the MSEs and the MSE matrices of shrinkage estimators improving on
the UMVUEs.  
In Section 5.1 we give some Monte Carlo studies to confirm the
theoretical results on the dominance relationship between the proposed 
estimators and the UMVUEs. 
In Section 5.2, we apply the proposed estimators of MSE matrices to
form confidence sets centered at shrinkage estimators and show their
usefulness through Monte Carlo studies.
We end this paper with some concluding remarks in Section 6.

\section{Preliminary facts on the estimation of MSE and MSE matrix for
 shrinkage estimators}
\label{sec:2}
Let $\bm{X}$ and $S$ be random variables which are independently
distributed as 
\begin{equation}
 \label{model}
 \bm{X} \sim \mathrm{N}_p(\bm{\theta}, \sigma^2 \bm{I}_p), \quad 
  S \sim \sigma^2\chi_n^2, 
\end{equation}
where $\mathrm{N}_p(\bm{\theta}, \sigma^2 \bm{I}_p)$ denotes 
$p$-variate normal distribution with a mean $\bm{\theta}$ and 
a covariance matrix $\sigma^2\bm{I}_p$ and $\chi_n^2$ denotes a
chi-square variable with $n$ degrees of freedom. 
Suppose that both $\bm{\theta}$ and $\sigma^2$ are unknown and that 
$p \ge 3$.
We write $W = \Vert \bm{X} \Vert^2/S$.
In this article we consider the following class of shrinkage estimators
\begin{equation}
 \label{eq:1}
  \bdeltap 
  =
  \left(
  1-\frac{\phi(W)}{W}
  \right)\bm{X}
\end{equation}
improving on $\bm{X}$ under squared error loss.
When $\phi(W)= (p-2)/(n+2)$, $\bdeltap$ is the James-Stein estimator 
(James and Stein\cite{James-Stein})
$$
\bm{\delta}_{JS} 
=
\left(
1-\frac{p-2}{n+2}\cdot\frac{1}{W}
\right)\bm{X}.
$$
When $\phi(W)= \max(W,(p-2)/(n+2))$,  $\bdeltap$ is the positive-part Stein
estimator(Baranchik\cite{Baranchik}), 
$$
\bm{\delta}_{JS}^+
=
\max
\left(
0,\; 1-\frac{p-2}{n+2}\cdot\frac{1}{W}
\right)\bm{X}.
$$
We assume that $\phi(\cdot)$ and $\phi'(\cdot)$ are finite. 
Most shrinkage estimators in the class (\ref{eq:1}) which have
been proposed in literature satisfy the condition.  
The MSE matrix of $\bdeltap$ is defined by 
$
\bm{M}(\bdeltap)=
\mathrm{E}[(\bdeltap-\bm{\theta})(\bdeltap-\bm{\theta})']
$
and the MSE is its trace, $R(\bdeltap)=\mathrm{tr}\bm{M}(\bdeltap)$.
We note that since $\bdeltap$ improves on $\bm{X}$, 
$R(\bdeltap) \le p\sigma^2$.
In this article we consider the estimation of 
$R(\bdeltap)$ and $\bm{M}(\bdeltap)$. 
As Carter et al.\cite{Carter-etal} and Wan et al.\cite{Wan-etal} suggested, 
the estimated MSE and MSE matrix are useful as precision measures for
$\bdeltap$. 
Denote by $\hat{\bm{M}}_0(\bdeltap)$ the UMVUE of $\bm{M}(\bdeltap)$.
They proposed to use $\hat{\bm{M}}_0(\bdeltap)$ in a quadratic form
which defines a confidence set for $\bm{\theta}$, 
\begin{equation}
 \label{ellipsoid}
  Q(\bdeltap) = 
  (\bdeltap - \bm{\theta})'
  \hat{\bm{M}}_0^{-1}(\bdeltap)
  (\bdeltap - \bm{\theta}).
\end{equation}

Let 
$\bm{\Gamma}$ be the orthogonal matrix such that 
$\bm{\Gamma}\bm{X} = (\Vert \bm{X} \Vert,0,\ldots,0)'$
and 
let
$\bm{E}_{11}$ be the $p \times p$ matrix which has one for (1,1) element
and zero elsewhere. 
Wan et al.\cite{Wan-etal} showed that 
when both $\phi(\cdot)$ and $\phi'(\cdot)$ are continuous, 
$\hat{\bm{M}}_0(\bdeltap)$ is expressed by 
\begin{align}
 \label{UMVUE:matrix}
 \hat{M}_0(\bm{\delta}_\phi) &= 
 \frac{S}{n} \cdot \bm{I}_p -S \cdot g_1(W) \cdot \bm{I}_p
 +S \cdot g_2(W) \cdot \frac{\bm{X}\bm{X}'}{\Vert \bm{X} \Vert^2}
 +S \cdot \frac{\phi^2(W)}{W} 
 \cdot \frac{\bm{X}\bm{X}'}{\Vert \bm{X} \Vert^2}, \notag\\
 &= S \cdot \bm{\Gamma}'
 \left\{
 \frac{\bm{I}_p}{n} - g_1(W) \cdot \bm{I}_p
 +
 g_3(W) \cdot \bm{E}_{11}
 \right\}
 \bm{\Gamma}, 
\end{align}
where 
$$
g_i(W)
=W^{\frac{n}{2}}
  \left\{
  \int_W^{\infty}
  t^{-\frac{n}{2}} \cdot
  \frac{h_i(t)}{t}\mbox{d}t  
  \right\},
\mbox{\ \ }i=1,2,\quad 
g_3(W)=g_2(W) + \frac{\phi^2(W)}{W},
$$
$$
h_1(W) = \frac{\phi(W)}{W},\quad 
h_2(W) = 2\left(
\frac{\phi(W)}{W}-\phi'(W)
\right).
$$
The UMVUE of $R(\bdeltap)$ is obtained by 
\begin{align}
 \label{UMVUE:MSE}
 \hat{R}_0(\bdeltap) &= 
 \mathrm{tr}\hat{\bm{M}}_0(\bm{\delta}_\phi) 
 =\frac{pS}{n} - S \cdot g(W)
  +\frac{S\phi^2(W)}{W}, 
\end{align}
where 
$$
g(W) = W^{\frac{n}{2}}
\left\{
\int_W^{\infty}
t^{-\frac{n}{2}} \cdot
\frac{h(t)}{t}\mbox{d}t  
\right\}, \quad
h(W) \equiv (p-2)\frac{\phi(W)}{W} +2\phi'(W).
$$
In the case of $\bm{\delta}_{JS}$, 
Efron and Morris\cite{Efron-Morris} and Carter et al.\cite{Carter-etal} 
showed that 
$\hat{R}_0({\bm{\delta}}_{JS})$ and 
$\hat{\bm{M}}_0({\bm{\delta}}_{JS})$
are expressed by  
\begin{equation}
 \label{eq:MSE_JS}
  \hat{R}_0({\bm{\delta}}_{JS}) =
  \frac{pS}{n} -
  \left(
    \frac{p-2}{n+2}
  \right)^2    
  \frac{S}{W}, 
\end{equation}
\begin{equation}
 \label{eq:matrix_JS}
 \hat{M}_0({\bm{\delta}}^{JS})
 =
 S \cdot \bm{\Gamma}'
 \left[
 \frac{\bm{I}_p}{n} -
 \frac{p-2}{(n+2)^2}\frac{1}{W}
 \left(
 2\bm{I}_p - (p+2)\bm{E}_{11}
 \right)
 \right]  
 \bm{\Gamma}, 
\end{equation}
respectively.
We can see that (\ref{eq:MSE_JS}) and (\ref{eq:matrix_JS}) 
coincide with (\ref{UMVUE:MSE}) and (\ref{UMVUE:matrix}) 
with $\phi(W) = (p-2)/(n+2)$.
$g_1(\cdot)$ and $g_3(\cdot)$ of $\hat{\bm{M}}_0(\bdelta_{JS})$ is
written by  
\begin{equation}
 \label{g:JS}
  g_1(W) = \frac{2(p-2)}{(n+2)^2}\cdot
  \frac{1}{W}, \quad
  g_3(W) = \frac{(p-2)(p+2)}{(n+2)^2}\cdot
  \frac{1}{W}.
\end{equation}
In the case of $\bm{\delta}_{JS}^+$, $\phi'(\cdot)$ is not continuous.
In order to apply the argument in Wan et al.\cite{Wan-etal} to 
the derivation of $\hat{R}_0(\bdelta_{JS}^+)$ and
$\hat{\bm{M}}_0(\bm{\delta}_{JS}^+)$, it requires some modifications. 
$\hat{R}_0({\bm{\delta}}_{JS}^+)$ and 
$\hat{\bm{M}}_0({\bm{\delta}}_{JS}^+)$
are expressed as follows, 
\begin{equation}
 \label{eq:MSE_positive}
  \hat{R}_0({\bm{\delta}}_{JS}^+) =
  \left\{
   \begin{array}{ll}
    \displaystyle{
     -\frac{pS}{n} + SW + C_0 S W^{n/2}
     }, 
     & 
     \displaystyle{
     W \le \frac{p-2}{n+2}},\\
    \\
    \displaystyle{
     \hat{R}_0({\bm{\delta}}_{JS}),
     }
     & \text{otherwise}, 
   \end{array}
  \right.
\end{equation}
\begin{equation}
 \label{eq:matrix_positive}
 \hat{\bm{M}}_0({\bm{\delta}}_{JS}^+)
 =
  \left\{
   \begin{array}{ll}
    \displaystyle{
     S \cdot \bm{\Gamma}'
     \biggl[
     \frac{\bm{I}_p}{n} -
     \left(
      \frac{2}{n}
      - C_1 W^{n/2}
     \right)
     \bm{I}_p}
      \\
     \qquad\qquad\qquad +
      \displaystyle{
      \left(
       W + C_2 W^{n/2}
       \right)
      \bm{E}_{11} 
     \biggr]
     \bm{\Gamma}, } & 
     \displaystyle{
     W \le \frac{p-2}{n+2}},\\     
    \\
    \displaystyle{
     \hat{\bm{M}}_0({\bm{\delta}}_{JS}),
     }
     & \text{otherwise}, 
   \end{array}
  \right.
\end{equation}
where 
$$
C_0 = 
2
\left(
\frac{p}{n} - \frac{p-2}{n+2}
\right)
\cdot
\left(
\frac{p-2}{n+2}
\right)^{-n/2}, 
$$
$$
C_1 = 
2
\left(
\frac{1}{n} - \frac{1}{n+2}
\right)
\cdot
\left(
\frac{p-2}{n+2}
\right)^{-n/2}, \quad 
C_2 =
\frac{4}{n+2}
\left(
\frac{p-2}{n+2}
\right)^{-n/2}.
$$
The proofs of (\ref{eq:MSE_positive}) and (\ref{eq:matrix_positive})
are given in the Appendix.
$g_1(W)$ and $g_3(W)$ of $\hat{\bm{M}}_0(\bdelta_{JS}^+)$ for 
$W \le (p-2)/(n+2)$  
is written by 
\begin{equation}
 \label{g:JS^+}
  g_1(W) = \frac{2}{n} - C_1 W^{n/2},\quad
  g_3(W) = W + C_2 W^{n/2}.
\end{equation}
When $W > (p-2)/(n+2)$, 
$g_1(W)$ and $g_3(W)$ of $\hat{\bm{M}}_0(\bdelta_{JS}^+)$ is identical
to (\ref{g:JS}).

It is easy to see from (\ref{eq:MSE_JS}) and (\ref{eq:MSE_positive})
that both $\hat{R}_0(\bm{\delta}_{JS})$ and
$\hat{R}_0(\bm{\delta}_{JS}^+)$ 
have an undesirable property of taking negative values for small $W$.
When we use an estimated MSE matrix to form a confidence set as in 
(\ref{ellipsoid}), it should be positive definite.
However we can also see from (\ref{eq:matrix_JS}) and
(\ref{eq:matrix_positive}) that both $\hat{\bm{M}}_0(\bm{\delta}_{JS})$
and $\hat{\bm{M}}_0(\bm{\delta}_{JS}^+)$ are not positive definite
for small $W$. 
In this article we consider the estimation of $R(\bdeltap)$ and 
$\bm{M}(\bdeltap)$ from a decision theoretical viewpoint. 
Kubokawa and Srivastava\cite{Kubokawa-Srivastava} considered the
estimation of the risk reduction and the risk reduction matrix of 
$\bdelta_{JS}$
$$
R^*(\bdelta_{JS}) = p\sigma^2 - R(\bdelta_{JS}),\quad 
\bm{M}^*(\bdelta_{JS}) = \sigma^2\bm{I}_p - \bm{M}(\bdelta_{JS}).
$$
They provided a class of estimators $\hat{R}^*(\bdelta_{JS})$ and 
$\hat{\bm{M}}^*(\bdelta_{JS})$ 
improving on the UMVUEs of $R^*(\bdeltap)$ and $\bm{M}^*(\bdeltap)$, 
$$
\hat{R}^*_0(\bdelta_{JS}) = \frac{pS}{n} - \hat{R}_0(\bdeltap),\quad
\hat{\bm{M}}^*_0(\bdelta_{JS}) = \frac{S}{n}\bm{I}_p - \hat{\bm{M}}_0(\bdeltap)
$$
under quadratic losses and satisfying 
$\hat{R}(\bdelta_{JS}) = (pS/n)-\hat{R}^*(\bdelta_{JS}) > 0$ and 
$\hat{\bm{M}}(\bdelta_{JS})=(S/n)\bm{I}_p - \hat{\bm{M}}^*(\bdelta_{JS})
> \bm{0}$. 
In the following sections we extend the argument in Kubokawa and
Srivastava\cite{Kubokawa-Srivastava} to a wider class of shrinkage
estimators $\bdeltap$ and propose some classes of nonnegative and
positive estimators of $R(\bdeltap)$ and nonnegative definite and
positive definite estimators of $\bm{M}(\bdeltap)$ improving on the
UMVUEs quadratic loss functions. 

\begin{remark}
 Carter et al.\cite{Carter-etal} and Wan et al.\cite{Wan-etal}
 considered this problem in the context of linear regression models, 
 $$
 \bm{Y} = \bm{A}\bm{\beta} + \bm{\epsilon}, 
 $$
 where $\bm{Y}$ is the $N$ observations, $\bm{A}$ is the $N \times p$
 design matrix and 
 $\bm{\epsilon} \sim \mathrm{N}_N(\bm{0}, \sigma^2\bm{I}_N)$ is the
 disturbance.  As mentioned in Kubokawa and
 Srivastava\cite{Kubokawa-Srivastava},  the setting (\ref{model})
 considered here can  be interpreted as a canonical form of linear
 regression models. 
 Let $\bm{B} = (\bm{A}'\bm{A})^{1/2}$ denote a $p \times p$
 matrix such that $\bm{B}\bm{B}'=\bm{A}'\bm{A}$.
 Then $\bm{X}$, $\bm{\theta}$, $S$ and $n$ corresponds to 
 $\bm{B}^{-1}\bm{A}'\bm{Y}$, $\bm{B}\bm{\beta}$, 
 $\bm{Y}(\bm{I}_N - \bm{A}(\bm{A}'\bm{A})^{-1}\bm{A}')\bm{Y}$ and $N-p$,
 respectively.
\end{remark}

\section{Improved estimators of the MSE}
\label{sec:3}
\subsection{An improved nonnegative estimator of the MSE}
\label{sec:3.1}
In this section we provide a nonnegative estimator of
$R(\bm{\delta}_{\bm{\phi}})$ improving on 
$\hat{R}_0(\bm{\delta}_{\bm{\phi}})$.
For evaluating an estimator $\hat{R}(\bdeltap)$, 
we use the following quadratic loss function, 
\begin{equation}
 \label{loss:MSE}
  L(\hat{R}(\bdeltap):R(\bdeltap))=
  (\hat{R}(\bdeltap) - R(\bdeltap))^2.
\end{equation}
Under the loss (\ref{loss:MSE}), the truncated estimator 
$\hat{R}^{TR}(\bdeltap)=
\max(0,\hat{R}_0(\bdeltap))$   
obviously dominates $\hat{R}_0(\bdeltap)$. 
The estimator proposed in this section is shown to dominates 
$\hat{R}^{TR}(\bdeltap)$

Let $a(W)$ be defined by 
$$
a(W)=\frac{n}{p}
\left(
g(W)-\frac{\phi^2(W)}{W}
\right).
$$
Then $\hat{R}_0(\bm{\delta}_{\bm{\phi}})$ is rewritten by 
\begin{equation}
 \label{a(W)}
  \hat{R}_0(\bm{\delta}_{\bm{\phi}})=
  pS
  \left(
   1-a(W)
  \right)/n.
\end{equation}
We introduce the following class of estimators with a function 
$\psi(\cdot)$, 
\begin{equation}
 \label{class:MSE}
  \hat{R}(\psi;\bm{\delta}_{\bm{\phi}})=
  \frac{pS}{n}
  \left(
   1 - a(W)\psi(W)
  \right),\quad 
  \psi(W) \le \frac{1}{a(W)}.
\end{equation}
Then we obtain the following theorem.
\begin{theorem}
 \label{th:MSE}
 Let $\psi_0(W)$ be defined by 
 $$
 \psi_0(W) = 
 \max\left[
 \min\left(
 1, 
 \frac{1}{a(W)}
 \right),\;
 \frac{1}{a(W)}
 \left(
 1-
 \frac{n(1+W)}{n+p+2}
 \right)
 \right].
 $$
 Then 
 \begin{align*}
  \hat{R}(\psi_0;\bm{\delta}_{\bm{\phi}})
  &=
  \frac{pS}{n}
  \left(
  1 - a(W)\psi_0(W)
  \right)\\
  &=
  \min\left[
  \max(\hat{R}_0(\bm{\delta}_{\bm{\phi}}), 0),\;
  \frac{pS(1+W)}{n+p+2}
  \right]
 \end{align*}
 improves on \(\hat{R}_0(\bm{\delta}_{\bm{\phi}})\) under the loss function 
 (\ref{loss:MSE}). 
\end{theorem}

\begin{proof}
 In the proof of this theorem, we use the similar procedure which
 Stein\cite{Stein-1964} used to derive the improved truncated
 estimator of 
 $\sigma^2$.  
 From the definition of 
 $\hat{R}(\psi;\bm{\delta}_{\bm{\phi}})$, 
 we have 
 \begin{align}
  \label{convex_psi}
  & (\hat{R}(\psi;\bm{\delta}_{\bm{\phi}}) - 
  R(\bm{\delta}_{\bm{\phi}}))^2\notag\\
  & \quad =
  \left(
  \frac{pS}{n}\cdot a(W) 
  \right)^2 \psi^2(W)
  -
  2 \left\{
  \left(
  \frac{pS}{n}
  \right)^2 a(W) 
  -
  \frac{pS}{n} \cdot R(\bm{\delta}_{\bm{\phi}})
  \right\}
  \psi(W)\notag\\
  &\qquad + 
  \left(
  \frac{pS}{n}-
  R(\bm{\delta}_{\bm{\phi}})
  \right)^2
 \end{align}
 Let $\psi^*(W)$ be the function which minimizes
 $
 \mathrm{E}
 [(\hat{R}(\psi;\bm{\delta}_{\bm{\phi}}) - 
 R(\bm{\delta}_{\bm{\phi}}))^2 \mid W]
 $.  
 From (\ref{convex_psi}), $\psi^*(W)$ is expressed by 
 \begin{align}
  \label{psi^*}
  \psi^*(W)
  &=
  \frac{1}{a(W)}
  -
  \frac{n\mathrm{E}[S/\sigma^2|W]}{\mathrm{E}[(S/\sigma^2)^2|W]}
  \cdot \frac{1}{a(W)}
   \cdot
   \frac{R(\bm{\delta}_{\bm{\phi}})}{p\sigma^2}.
 \end{align}
 Denote by 
 $f_k(\cdot)$ and $f_k(\cdot ; \lambda)$ 
 the density functions of the central 
 and the non-central $\chi^2$ distributions with $k$ degrees of freedom
 and the non-centrality parameter 
 $\lambda=\Vert \bm{\theta} \Vert^2/\sigma^2$.
 Then we have 
 \begin{align}
  \label{ineq:chi^2}
  \frac{\mathrm{E}[S/\sigma^2|W]}{\mathrm{E}[(S/\sigma^2)^2|W]}
   &=
  \frac{\int u^2 f_p(uW;\lambda)f_n(u) du}
  {\int u^3 f_{p}(uW;\lambda)f_n(u) du}\notag\\
  &\leq
  \frac{\int u^2 f_p(uW)f_n(u) du}
   {\int u^3 f_p(uW)f_n(u) du} =
   \frac{1+w}{n+p+2}, 
 \end{align}
 from the monotone nondecreasingness of 
 $f_p(uW;\lambda)/f_p(uW)$ in $u$. 
 $R(\bdeltap) \le p\sigma^2$.
 Hence from (\ref{psi^*}), (\ref{ineq:chi^2}) and the fact that 
 $0 < R(\bdeltap) \le p\sigma^2$, 
 $\psi^*(W)$ satisfies 
 $$
 \frac{1}{a(W)}
 \left(
 1-
 \frac{n(1+W)}{n+p+2}
 \right)
 \leq
 \psi^*(W)
 \leq
 \frac{1}{a(W)}
 $$
 for any $\lambda$ and 
 any $W$ such that $a(W) > 0$.
 Denote
 ${\cal R^+}=\{W \mid W \geq 0\}$
 and define ${\cal A}$, ${\cal B}$ and ${\cal C}$ as follows,
$$
  \begin{array}{l}
   {\cal A} =
    \displaystyle{
    \left\{ 
     W \left| \frac{1}{a(W)} \leq 1, \; a(W)>0, \; W \ge 0
       \right.
    \right\}}\\
   \\
  {\cal B}
  =
  \displaystyle{
  \left\{ 
   W  \left| \frac{1}{a(W)}
       \left(
	1-
	\frac{n(1+W)}{n+p+2}
       \right) 
  \geq 1, \; a(W)>0, \; W \ge 0
  \right.
  \right\}}\\
   \\
  {\cal C}
  =
  {\cal R^+} \setminus ({\cal A} \cup {\cal B}).
 \end{array}
$$
 We note that 
 $\mathrm{E}[(\hat{R}(\psi;\bm{\delta}_{\bm{\phi}}) - 
 R(\bm{\delta}_{\bm{\phi}}))^2
 \mid W]$ 
 is convex on $\psi(W)$.
 Therefore, if we set 
 $$
 \psi(W)
 =\min
 \left(
 1,\frac{1}{a(W)}
 \right) \quad \text{for} \quad W \in {\cal A},
 $$
 $\psi^*(W)$ satisfies $\psi^*(W) \leq \psi(W) \leq 1$ for 
 $W \in {\cal A}$, which implies 
 \begin{equation}
  \label{cond:A}
 \mathrm{E}[
 (\hat{R}(\psi;\bm{\delta}_{\bm{\phi}})-R(\bm{\delta}_{\bm{\phi}}))^2
 \mid W \in {\cal A}] 
  \leq
  \mathrm{E}[(\hat{R}_0(\bm{\delta}_{\bm{\phi}})-R(\bm{\delta}_{\bm{\phi}}))^2|W \in {\cal A}].
 \end{equation}
 Similarly, if we set
 $$
 \psi(W)
 =
 \max\left[
 \frac{1}{a(W)}
 \left(
 1-
 \frac{n(1+W)}{n+p+2}
 \right),\; 1
 \right]  \quad \text{for} \quad W \in {\cal B},
 $$
 $\psi^*(W)$ satisfies $\psi^*(W) \geq \psi(W) \geq 1 $
 for $W \in {\cal B}$ and then 
 \begin{equation}
  \label{cond:B} 
   \mathrm{E}[
   (\hat{R}(\psi;\bm{\delta}_{\bm{\phi}})-R(\bm{\delta}_{\bm{\phi}}))^2
   \mid W  
   \in {\cal B}] 
 \leq
 \mathrm{E}[(\hat{R}_0(\bm{\delta}_{\bm{\phi}})-R(\bm{\delta}_{\bm{\phi}}))^2|W \in {\cal B}].
 \end{equation}
 For $W \in {\cal C}$, if we set \(\psi(W)=1\), 
 \begin{equation}
  \label{cond:C}
 \mathrm{E}[
 (\hat{R}(\psi;\bm{\delta}_{\bm{\phi}})-R(\bm{\delta}_{\bm{\phi}}))^2
 \mid W \in {\cal C}] 
 =
 \mathrm{E}[(\hat{R}_0(\bm{\delta}_{\bm{\phi}})-R(\bm{\delta}_{\bm{\phi}}))^2|W \in {\cal C}].
 \end{equation}
 By combining (\ref{cond:A}), (\ref{cond:B}) and (\ref{cond:C}), 
 we can complete the proof.
\end{proof}

We have shown that 
$\hat{R}(\psi_0;\bm{\delta}_{\bm{\phi}})$
dominates 
$\hat{R}_0(\bm{\delta}_{\bm{\phi}})$ under the loss (\ref{loss:MSE}).
We note that the unbiasedness of $\hat{R}_0(\bm{\delta}_{\bm{\phi}})$
was not used in the proof. 
By following the proof, we can see that any estimator in the class
(\ref{class:MSE}) which does not satisfy 
\begin{equation}
 \label{nec-cond}
 0 \leq \hat{R}(\psi;\bm{\delta}_{\bm{\phi}}) \leq
 \frac{pS(1+W)}{n+p+2} 
\end{equation}
is improved by the one which is truncated it to satisfy
(\ref{nec-cond}). 
In other words, (\ref{nec-cond}) is a necessary condition that 
$\hat{R}(\psi;\bm{\delta}_{\bm{\phi}})$ is admissible in the class 
(\ref{class:MSE}). 

In the case of $\bm{\delta}_{JS}$, 
$a(W)=n(p-2)^2/(p(n+2)^2 W)$ from (\ref{eq:MSE_JS}).
Then ${\cal B}$ for $\bm{\delta}_{JS}$ is written by 
$$
{\cal B} =
  \left\{ 
   W  \left| 
   \left(
   \frac{n+2}{p-2}
   \right)^2 \cdot \frac{pW}{n} \cdot
       \left(
	1-
	\frac{n(1+W)}{n+p+2}
       \right) 
  \geq 1, \; W \ge 0
  \right.
  \right\}.
$$
Hence if $(p,n)$ satisfy
$$
\left(
\frac{p+2}{n}
\right)^2
-
\frac{4(n+p+2)}{p}
\left(
\frac{p-2}{n+2}
\right)^2
\leq 0, 
$$
${\cal B} = \emptyset$ and then 
$\hat{R}(\psi_0;\bm{\delta}_{JS})$ is identical to 
$\hat{R}^{TR}(\bm{\delta}_{JS})$. 
In general, when $n$ and $p$ are large, 
${\cal B} = \emptyset$ and hence 
$\hat{R}(\psi_0;\bdeltap)$ is identical to 
$\hat{R}^{TR}(\bdeltap)$. 
However if ${\cal B} \neq \emptyset$,  
$\hat{R}(\psi_0;\bm{\delta}_{\bm{\phi}})$ is shown to 
improve on $\hat{R}^{TR}(\bm{\delta}_{\bm{\phi}})$. 

\begin{corollary}
 If $(p,n)$ satisfies ${\cal B} \neq \emptyset$, 
 $\hat{R}(\psi_0;\bm{\delta}_{\bm{\phi}})$
 improves on  $\hat{R}^{TR}(\bm{\delta}_{\bm{\phi}})$
 under the loss (\ref{loss:MSE}). 
\end{corollary}

\begin{proof}
 As mentioned above, $\hat{R}^{TR}(\bm{\delta}_{\bm{\phi}})$ is
 identical to $\hat{R}(\psi;\bm{\delta}_{\bm{\phi}})$ with 
 $$
 \psi(W)=
 \left\{
 \begin{array}{ll}
  \displaystyle{
   \min\left(
	1,\; \frac{1}{a(W)}
       \right), 
   }& \text{if } W \in {\cal A}, \\
  \\
  1, & \text{otherwise}.
 \end{array}
 \right.
 $$
 By following the proof of Theorem \ref{th:MSE}, we can see that 
 $\hat{R}(\psi_0;\bm{\delta}_{\bm{\phi}})$
 improves on  $\hat{R}^{TR}(\bm{\delta}_{\bm{\phi}})$
 under the loss (\ref{loss:MSE}) when ${\cal B} \neq \emptyset$. 
\end{proof}

We present in Figure 1 the shape of 
$1-a(W)\psi_0(W)$ for $\delta_{JS}$ and $\delta_{JS}^+$
with $(p,n)=(5,5)$.
The dotted line and the solid line represent $1-a(W)$ and
$1-a(W)\psi_0(W)$, respectively. 
We can see that 
when $\hat{R}_0(\bm{\delta}_{\bm{\phi}})<0$, 
$\hat{R}(\psi_0;\bm{\delta}_{\bm{\phi}})$
shrinks towards 0 and even when 
$\hat{R}_0(\bm{\delta}_{\bm{\phi}}) > 0$, 
$\hat{R}(\psi_0;\bm{\delta}_{\bm{\phi}})$ 
slightly shrinks it towards 0.

\begin{figure}[htbp]
 \centering
 \begin{tabular}{cc}
  \includegraphics[scale=0.7]{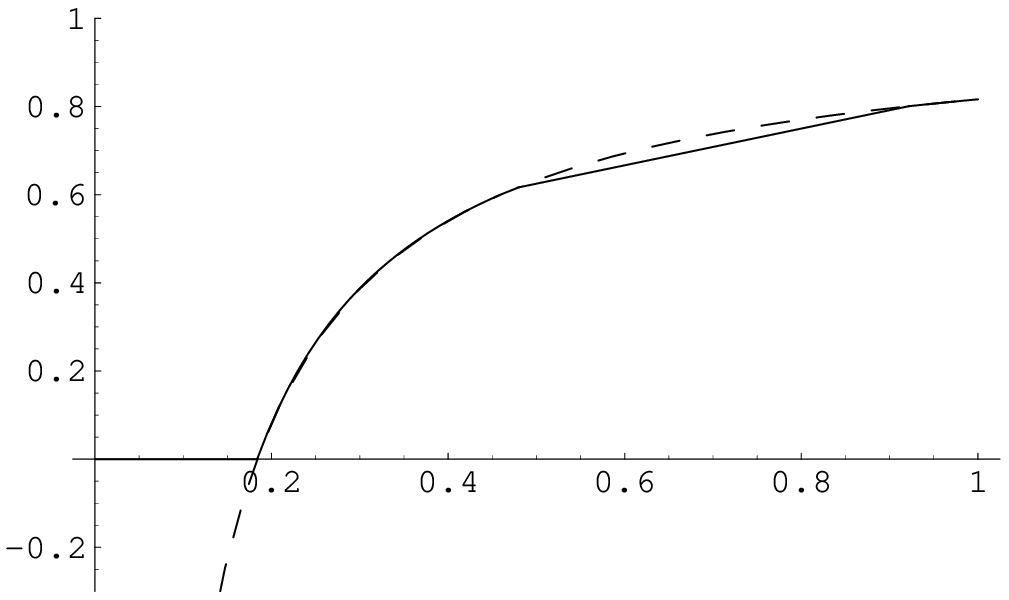} &
  \includegraphics[scale=0.7]{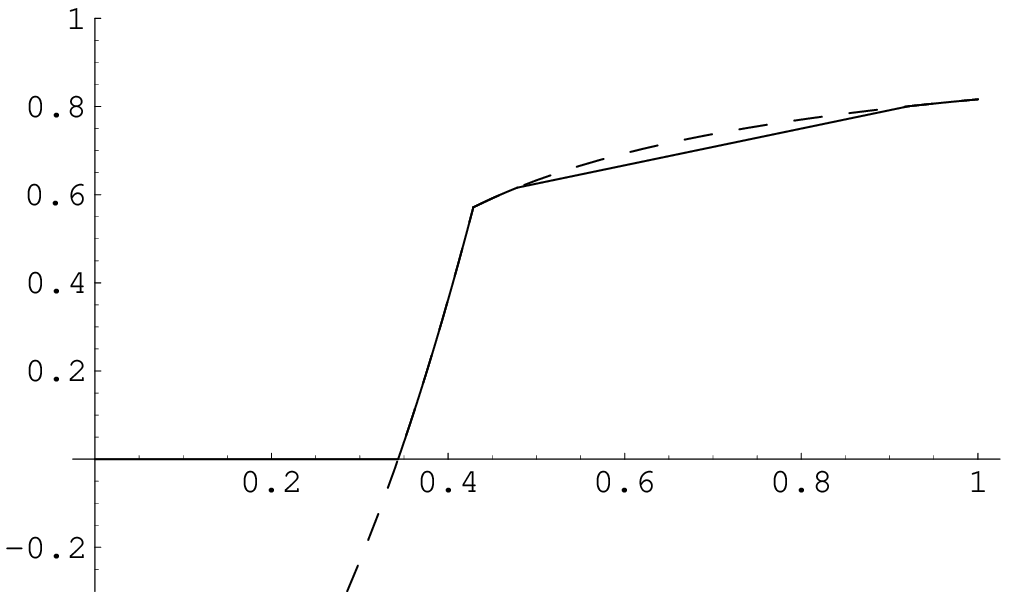}\\
  $1-a(W)\psi_0(W)$ for $\delta_{JS}$ & 
  $1-a(W)\psi_0(W)$ for $\delta_{JS}^+$
  \end{tabular}
 \caption{\(1-a(W)\psi_0(W)\) of 
 \(\bm{\delta}_{JS}\) and \(\bm{\delta}_{JS}^+\) when $(p,n)=(5,5)$}
 \label{fig:1}
\end{figure}

\subsection{Improved estimation of the risk reduction and positive 
  estimation of MSE} 
\label{sec:3.2}
In the previous section we have shown that
$\hat{R}(\psi_0;\bm{\delta}_{\bm{\phi}})$ improves 
on $\hat{R}_0(\bm{\delta}_{\bm{\phi}})$. 
However we can easily see that 
$\hat{R}(\psi_0;\bm{\delta}_{\bm{\phi}})$ can take zero with 
positive probability. From a practical viewpoint, 
$\hat{R}(\psi_0;\bm{\delta}_{\bm{\phi}})$ may still be
undesirable. 
In order to obtain positive estimators, 
we take an similar approach to the one in Kubokawa and
Srivastava\cite{Kubokawa-Srivastava} and consider 
the improved estimation of the risk reduction 
$R^*(\bdeltap) = p\sigma^2 - R(\bdeltap)$.
In accordance with Kubokawa and Srivastava\cite{Kubokawa-Srivastava},  
for evaluating the estimator $\hat{R}^*(\bdeltap)$, 
we use the following quadratic loss function,  
\begin{equation}
 \label{loss:reduction-MSE}
  L(\hat{R}^*(\bdeltap) ; R^*(\bdeltap))
  =
  (\hat{R}^*(\bdeltap) - R^*(\bdeltap))^2.
\end{equation}
From (\ref{a(W)}), the UMVUE of $R^*(\bdeltap)$ is 
$\hat{R}^*_0(\bdeltap) = pSa(W)/n$.
In order to improve on $\hat{R}^*_0(\bdeltap)$, 
we introduce the class of estimators 
$\hat{R}^*(\psi;\bdeltap)=pSa(W)\psi(W)/n$.
We can see that 
$pS/n - \hat{R}^*(\psi;\bdeltap)$ is identical to 
$\hat{R}(\psi;\bdeltap)$ defined in (\ref{class:MSE}).
It is well known that $R^*(\bdeltap)/\sigma^2$  depends only on
$\lambda$(e.g. Efron and Morris\cite{Efron-Morris}). 
We suppose that 
$R^*(\bdeltap)/\sigma^2$ is maximized at $\lambda=0$.
If $\phi(W)/W$ is nonincreasing, $\bdeltap$ satisfy this
condition(e.g. Kubokawa\cite{Kubokawa-1988}, Casella\cite{Casella}).
Many shrinkage estimators including $\bdelta_{JS}$ and $\bdelta_{JS}^+$
satisfy this condition. 
Denote $R^*(\bdeltap)/\sigma^2$ with $\lambda=0$
by $\alpha_{p,n}(\bdeltap)$.
By using the similar argument in Theorem 1, 
we can obtain the following result. 

\begin{theorem}
 \label{th:MSE-2}
 Let $W_{p,n}(\bdeltap)$ be the solution of the equation
 \begin{equation}
  \label{eq:W}
   \frac{1+W}{a(W)} = \frac{p(n+p+2)}{n} \cdot
   \frac{1}{\alpha_{p,n}(\bdeltap)}.
 \end{equation}
 Define $\psi_1(W)$ by 
 $$
 \psi_1(W) = 
 \max\left(
 1,\; \frac{n}{p} \cdot \frac{1+W_{p,n}(\bdeltap)}{n+p+2}
 \cdot \frac{\alpha_{p,n}(\bdeltap)}{a(W)}
 \right).
 $$
 Then $\hat{R}^*(\psi_1;\bdeltap)$ improves on 
 $\hat{R}^*_0(\bdeltap)$.
\end{theorem}
The proof is similar to the proof of Theorem \ref{th:MSE} and is given
in the Appendix. 

By using the similar argument in Kubokawa and
Srivastava\cite{Kubokawa-Srivastava}, we can also obtain the following
result. 

\begin{theorem}
 \label{th:MSE-3}
  Suppose that 
 \begin{enumerate}
  \item[(i)] $a(W)$ is nonincreasing ;
  \item[(ii)] $\psi(W)$ is nondecreasing and  
	     $\lim_{W \to \infty} \psi(W) = 1$ ;
  \item[(iii)] $\psi(W) \ge \psi_2(W)$ where
	     $$
	     \psi_2(W) = 
	     \min\left(
	     1, \; 
	     \frac{\alpha_{p,n}(\bdeltap)}
	     {p a(W)} 
	     \right).
	     $$
 \end{enumerate}
 Then $\hat{R}^*(\psi; \bdeltap)$ improves on 
 $\hat{R}^*_0(\bdeltap)$ under the loss (\ref{loss:reduction-MSE}).   
 Especially, $\hat{R}^*(\psi_2; \bdeltap)$ improves on 
 $\hat{R}^*_0(\bdeltap)$.
\end{theorem}

The proof of this theorem is also given in the Appendix. 
We note that we can see from (\ref{eq:MSE_JS}) and
(\ref{eq:MSE_positive}) that 
$\bdelta_{JS}$ and $\bdelta_{JS}^+$ satisfies the condition (i).
The results of Theorem \ref{th:MSE-2} and Theorem \ref{th:MSE-3} lead us
to propose the following estimators of $\hat{R}(\bdeltap)$, 
\begin{align*}
 \hat{R}(\psi_1; \bdeltap)
 &=\frac{pS}{n} - 
 \hat{R}^*(\psi_1; \bdeltap)\\
 &= \max 
 \left[
 \frac{pS}{n}(1-a(W)), \;
 \frac{pS}{n}\left(
 1 - \frac{n(1+W_{p,n}(\bdeltap))}{n+p+2} \cdot 
 \frac{\alpha_{p,n}(\bdeltap)}{p}
 \right)
 \right], \\
 \hat{R}(\psi_2; \bdeltap)
 &=\frac{pS}{n} - 
 \hat{R}^*(\psi_2; \bdeltap)\\
 &= \max 
 \left[
 \frac{pS}{n}(1-a(W)), \;
 \frac{pS}{n}\left(
 1 - \frac{\alpha_{p,n}(\bdeltap)}{p}
 \right)
 \right].
\end{align*}
Since $\alpha_{p,n}(\bdeltap) < p$, 
$\hat{R}(\psi_2; \bdeltap)$ is always positive.
If $p$ and $n$ satisfy
\begin{equation}
 \label{ineq:W0}
  \gamma_{p,n}(\bdeltap) := 
  \frac{n(1+W_{p,n}(\bdeltap))}{n+p+2} <  \frac{p}{\alpha_{p,n}(\bdeltap)(\bdeltap)}, 
\end{equation}
$\hat{R}(\psi_1; \bdeltap) > 0$ with probability one.
We can prove 
$\gamma_{p,n}(\bdelta_{JS}) < p/\alpha_{p,n}(\bdeltap)(\bdelta_{JS})$
analytically. 
The proof is given in the Appendix.

We present $\gamma_{p,n}(\bdelta_{JS})$, $\gamma_{p,n}(\bdelta_{JS}^+)$,
$W_{p,n}(\bdelta_{JS})$ and $W_{p,n}(\bdelta_{JS}^+)$ for 
$(p,n) = (5,5)$, $(10,5)$, $(5,10)$ and $(10,10)$
in Table \ref{table:gamma} and \ref{table:w0-1}. 
We can see $\gamma_{p,n}(\bdelta_{JS}) < 1$ and 
$\gamma_{p,n}(\bdelta_{JS}^+) < 1$ for all cases.
Since $p/\alpha_{p,n}(\bdeltap)(\bdeltap) > 1$,
$\gamma_{p,n}(\bdelta_{JS})$ and $\gamma_{p,n}(\bdelta_{JS}^+)$
satisfy (\ref{ineq:W0}).
Many shrinkage estimators seem to satisfy (\ref{ineq:W0}). 
However it seems to be difficult to prove (\ref{ineq:W0}) 
for general $\bdeltap$ analytically.  
We note that $\gamma_{p,n}(\bdelta_{JS}) < 1$ and 
$\gamma_{p,n}(\bdelta_{JS}^+) < 1$ mean that $\psi_1(\cdot)$ does not
satisfy the conditions in Theorem \ref{th:MSE-3}.
Hence from Table \ref{table:gamma}, 
$\hat{R}(\psi_1;\bdeltap)$ is not necessarily included in the class
of estimators in Theorem \ref{th:MSE-3}.

\begin{table}[htbp]
 \centering
 \caption{$\gamma_{p,n}(\bdelta_{JS})$ and
 $\gamma_{p,n}(\bdelta_{JS}^+)$ for some $p$ and $n$} 
 \label{table:gamma}
 \begin{tabular}{ccccc}\hline
  $(p,n)$ & $(5,5)$ & $(10,5)$ & $(5,10)$ & $(10,10)$\\ \hline
  $\gamma_{p,n}(\bdelta_{JS})$  & 
  0.6795 & 0.7452 & 0.7774 & 0.8228\\
  $\gamma_{p,n}(\bdelta_{JS}^+)$ &
  0.6399 & 0.7056 & 0.7484 & 0.7921\\ \hline
 \end{tabular}
\end{table}
\begin{table}[htbp]
 \centering
 \caption{$W_{p,n}(\bdelta_{JS})$ and $W_{p,n}(\bdelta_{JS}^+)$ for some
 $p$ and $n$}  
 \label{table:w0-1}
 \begin{tabular}{ccccc}\hline
  $(p,n)$ & $(5,5)$ & $(10,5)$ & $(5,10)$ & $(10,10)$\\ \hline
  $W_{p,n}(\bdelta_{JS})$  & 
  0.6307 & 1.533 & 0.3216 & 0.8102\\
  $W_{p,n}(\bdelta_{JS}^+)$ &
  0.5357 & 1.399 & 0.2722 & 0.7427\\ \hline
 \end{tabular}
\end{table}

So far we considered the improved estimation of $R^*(\bdeltap)$. 
Next we consider to evaluate the resulting estimators
$\hat{R}(\psi_1; \bdeltap)$ and $\hat{R}(\psi_2; \bdeltap)$ 
as estimators of $R(\bdeltap)$. 
Based on the fact that (\ref{nec-cond}) is a necessary condition that  
$\hat{R}(\psi;\bm{\delta}_{\bm{\phi}})$ is admissible in the class 
(\ref{class:MSE}), 
we can provide estimators improving on 
$\hat{R}(\psi_1; \bdeltap)$ and $\hat{R}(\psi_2; \bdeltap)$
under the loss (\ref{loss:MSE}).

\begin{theorem}
 \label{th:MSE-5}
 Suppose that ${\cal B} \neq \emptyset$. 
 Let $\psi_1^{TR}(W)$ and $\psi_2^{TR}(W)$ 
 be defined by
 $$
 \psi_j^{TR}(W) = \max
 \left[
 \psi_j(W),\;
 \frac{1}{a(W)}
 \left(
 1 - \frac{n(1+W)}{n+p+2}
 \right)
 \right]. 
 $$
 Then 
 $\hat{R}(\psi_j^{TR}; \bdeltap)$ improves on $\hat{R}(\psi_j;
 \bdeltap)$ under the loss (\ref{loss:MSE}). 
\end{theorem}

Based on the proof of Theorem \ref{th:MSE}, the proof is easy and
omitted. 
It is clear that 
$\hat{R}(\psi_j^{TR} ; \bdeltap) > 0$ when 
$\hat{R}(\psi_j ; \bdeltap) > 0$.
The dominance relation between 
$\hat{R}(\psi_j^{TR} ; \bdeltap)$ and 
$\hat{R}_0(\bdeltap)$ is interesting.
But it seems difficult to clarify it at this point.

\section{Improved estimators of the MSE matrix}
\subsection{An improved nonnegative definite estimation of the MSE
  matrix} 
\label{sec:4.1}
In this section we provide a nonnegative definite estimator of
$\bm{M}(\bm{\delta}_\phi)$ improving on 
$\hat{\bm{M}}_0(\bm{\delta}_\phi)$.  
For evaluating an estimator $\hat{\bm{M}}(\bdeltap)$, 
we use the following squared loss functions, 
\begin{equation}
 \label{loss:matrix}
  L(\hat{\bm{M}}(\bdeltap)
  ;\bm{M}(\bdeltap))=
  \mathrm{tr}
  (\hat{\bm{M}}(\bdeltap) 
  - \bm{M}(\bdeltap))^2.
\end{equation}
We note that 
$L(\hat{\bm{M}}(\bdeltap) ; \bm{M}(\bdeltap))$
is the sum of quadratic losses of each element. 
In order to consider the improvement of
$\hat{\bm{M}}_0(\bm{\delta}_\phi)$, 
we introduce the following class of estimators with functions
$\xi(\cdot)$ and $\eta(\cdot)$, 
\begin{align}
 \label{class:matrix}
 \hat{\bm{M}}(\xi,\eta ; \bm{\delta}_{\bm{\phi}})
 &= S \cdot \bm{\Gamma}'
 \Biggl\{
 \frac{\bm{I}_p}{n} - g_1(W)\xi(W)
 \cdot 
 (\bm{I}_p - \bm{E}_{11})
 \notag\\
 &\qquad\qquad\qquad + 
 \bigl(
 g_3(W) - g_1(W)\eta(W) 
 \bigr)
 \cdot \bm{E}_{11}
 \Biggr\}
 \bm{\Gamma}\notag\\
 &=
 S \bm{\Gamma}'\bm{L}(\xi,\eta;\bdeltap)\bm{\Gamma}, 
\end{align}
where 
\begin{align*}
\bm{L}(\xi,\eta;\bdeltap) &= 
\mathrm{diag}
(l_0(\eta;\bdeltap),l_1(\xi;\bdeltap), 
\ldots, l_1(\xi;\bdeltap)), 
\end{align*}
$$
l_0(\eta;\bdeltap) = 
\frac1n - g_1(W)\eta(W) + g_3(W),\quad
l_1(\xi;\bdeltap) = 
\frac1n - g_1(W)\xi(W).
$$
It is clear that $\hat{\bm{M}}(\xi,\eta;\bm{\delta}_{\bm{\phi}})$ with
$\xi(W)=1$ and $\eta(W)=1$ is identical to
$\hat{\bm{M}}_0(\bm{\delta}_{\bm{\phi}})$. 
$\hat{\bm{M}}(\xi,\eta;\bm{\delta}_{\bm{\phi}})$ with
$\eta(W)=1$ coincides with the class which Kubokawa and
Srivastava\cite{Kubokawa-Srivastava} considered for the James-Stein
estimator.
In the case of positive-part Stein estimator, however, 
$g_3(W)-g_1(W)$ is written by 
\begin{align*}
 g_3(W) - g_1(W) &=
 \left\{
 \begin{array}{ll}
  \displaystyle{
   \left(
    \frac{2}{n+2} + \frac{2}{n}
   \right)
   W^{n/2}
   \left(
   \frac{p-2}{n+2}
   \right)^{-n/2} + W
   - \frac{2}{n}
   }, & 
   \displaystyle{
   W \le \frac{p-2}{n+2}, 
   }\\
  \\
  \displaystyle{
  \frac{p(p-2)}{(n+2)^2} \cdot \frac{1}{W}
  }, & \text{otherwise}.
 \end{array}
\right.
\end{align*}
Denote $l_0(\eta:\bdelta_{JS}^+)$ with $\eta(W)=1$ by
$l_0(\bdelta_{JS}^+)$.  
Then we can see that $l_0(\bdelta_{JS}^+) < 0$ for small $W$.
This means that any estimator $\hat{\bm{M}}(\xi,\eta ; \bdelta_{JS}^+)$ 
with $\eta(W) = 1$ does not always take nonnegative definite. 
So we consider to find improved estimators from the class 
(\ref{class:matrix}).

\begin{figure}[t]
 \centering
 \begin{tabular}{cc}
 \includegraphics[scale=0.35]{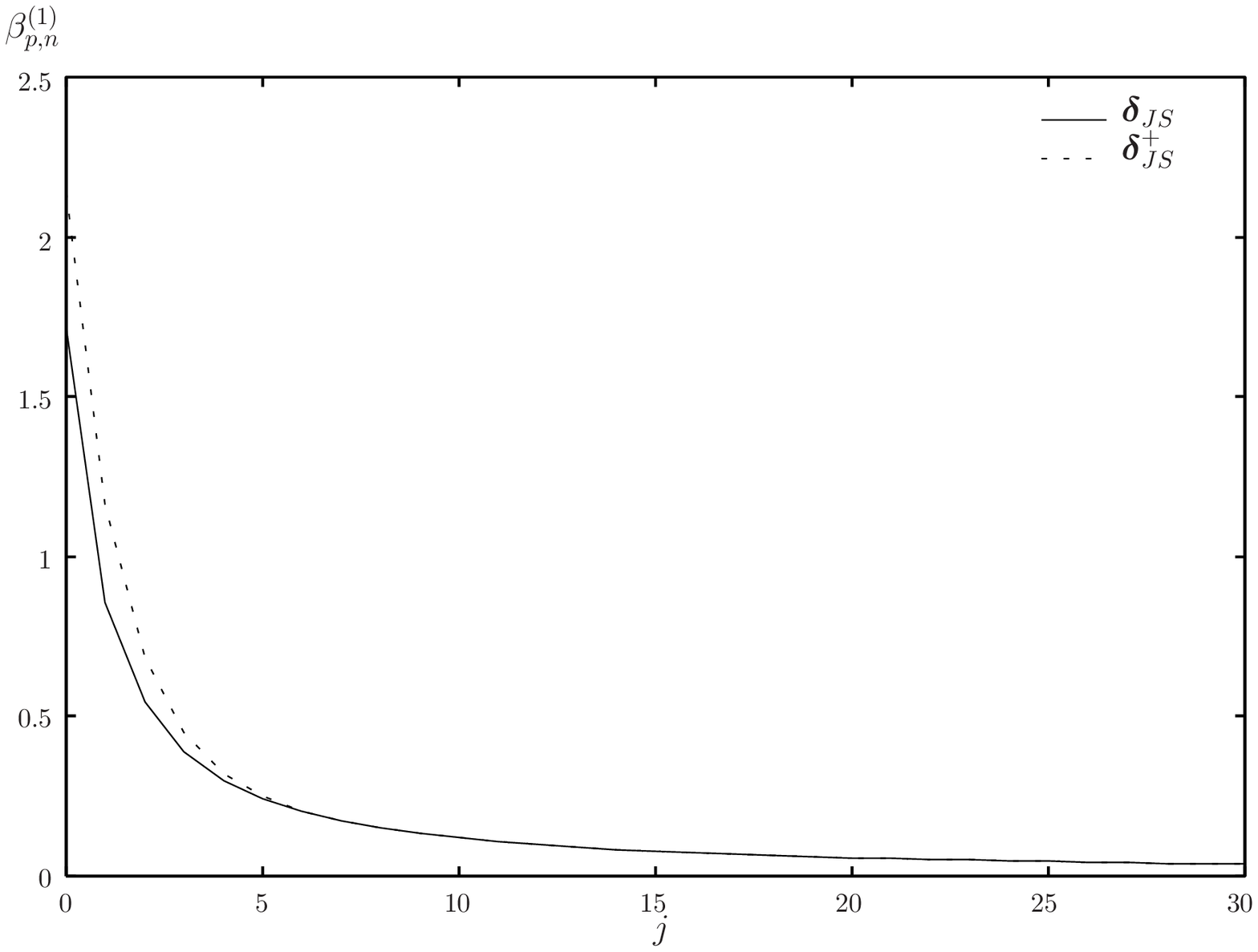} &
 \includegraphics[scale=0.35]{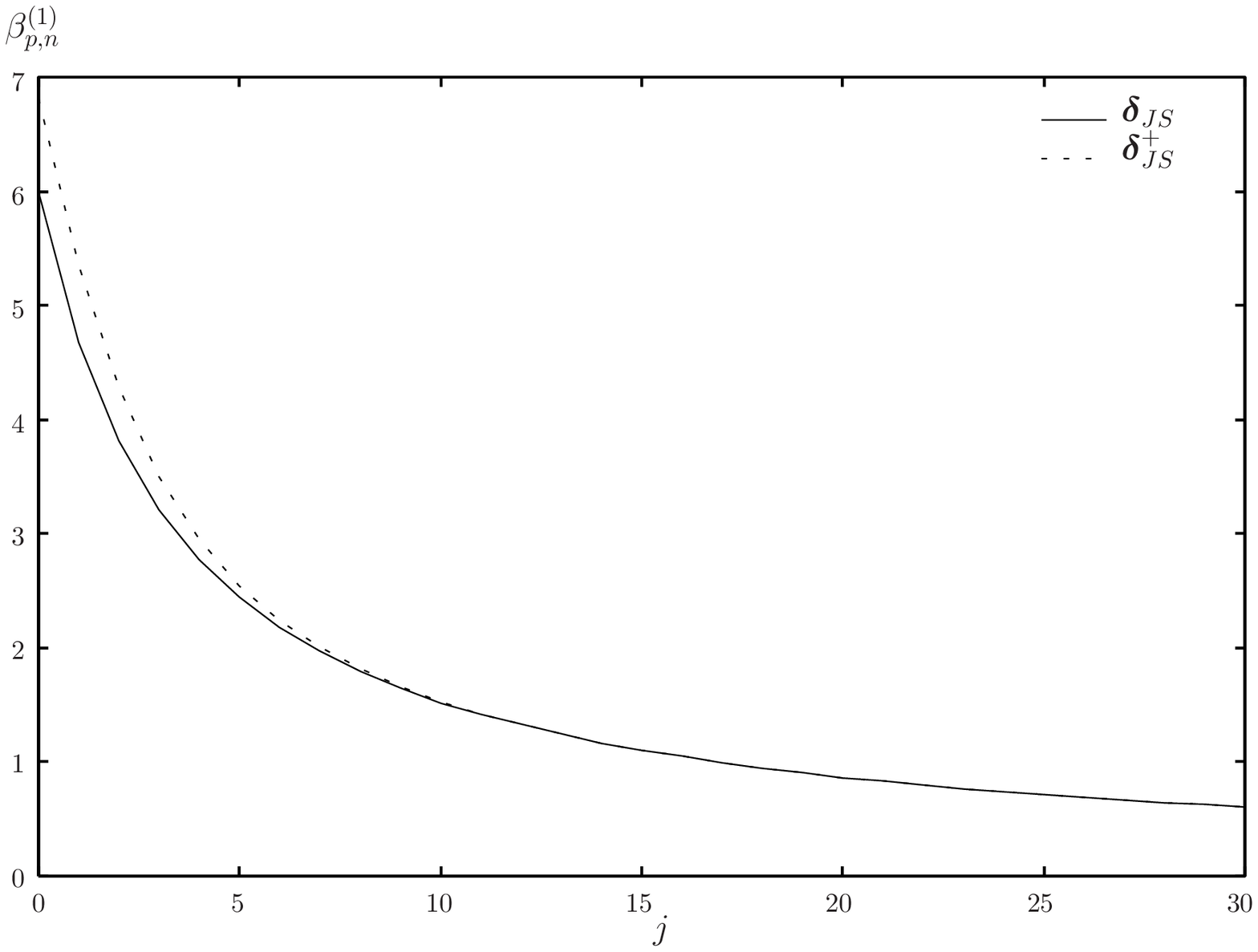}\\
 (i) $(p,n)=(5,5)$ & (ii) $(p,n)=(10,10)$ 
 \end{tabular}
 \caption{$\beta^{(1)}_{p,n}(j:\bdeltap)$ of $\bdelta_{JS}$ and
 $\bdelta_{JS}^+$ with $(p,n)=(5,5)$, $(10,10)$}
\label{figure:beta-1}
\end{figure}

We assume that $g_1(W)$ is nonincreasing.
From (\ref{g:JS}) and (\ref{g:JS^+}), we can see that 
$\bdelta_{JS}$ and $\bdelta_{JS}^+$ satisfy this condition.
Let $u$ and $v$ be random variables which distribute as chi-square
distribution with $p+2j$ and $n$ degrees of freedom, respectively.
Then define $\beta^{(1)}_{p,n}(\bdeltap)$ 
$$
\beta^{(1)}_{p,n}(\bdeltap) :=
\inf_{j \ge 0}\beta^{(1)}_{p,n}(j;\bdeltap), 
$$
where 
$$
\beta^{(1)}_{p,n}(j;\bdeltap) =
\mathrm{E}\left[
\frac{2(p-1)\phi(u/v)}{u/v} - 
\frac{(p+2j-1)b(u/v)}{p+2j}\cdot 
\right], 
$$
$$
  b(W) = \frac{4 \phi(W)}{W} + \frac{(n+2)\phi^2(W)}{W} - 4\phi'(W)
  - 4\phi(W)\phi'(W).
$$
Figure \ref{figure:beta-1} represents the behavior of 
$\beta^{(1)}_{p,n}(j;\bdeltap)$ of $\bdelta_{JS}$ and $\bdelta_{JS}^+$ 
with $(p,n)=(5,5)$, $(10,10)$ computed by Monte Carlo simulation. 
The result of Figure \ref{figure:beta-1} indicates that
$\beta^{(1)}_{p,n}(\bdelta_{JS}) \ge 0$ and 
$\beta^{(1)}_{p,n}(\bdelta_{JS}^+) \ge 0$.
We obtained the same results numerically with respect to other $(p,n)$
and other $\bdeltap$. 

By using the similar argument in the proof of Theorem \ref{th:MSE}, 
we can obtain the following theorem.

\begin{theorem}
 \label{th:matrix-1}
 Suppose that $\beta^{(1)}_{p,n}(\bdeltap) \ge 0$.
 Let $\xi_0(W)$ and $\eta_0(W)$ be defined by 
 \begin{align*}
 \xi_0(W) &= 
 \max\left[
 \min\left(
 1, \frac{1}{n g_1(W)}
 \right),\;
 \frac{1}{g_1(W)}
 \left(
 \frac{1}{n} - 
 \frac{1+W}{n+p+1}
 \right)
 \right], \\
 \eta_0(W) &= 
 \min\left[
 1, \; \frac{1}{g_1(W)}
 \left(
 \frac{1}{n} + g_3(W)
 \right)
  \right]
 \end{align*}
 Then 
 $\hat{\bm{M}}(\xi_0,\eta_0;\bdeltap)$
 improves on 
 $\hat{\bm{M}}_0(\bdeltap)$ under the loss (\ref{loss:matrix}). 
\end{theorem}

The proof is similar to the one of Theorem \ref{th:MSE}
and is given in the Appendix. 
By following the proof of the theorem and by using the same argument in 
Section \ref{sec:3.1}, we can see that 
\begin{equation}
 \label{ineq:xi^*}
  \frac{1}{g_1(W)}
  \left(
   \frac{1}{n} -
   \frac{1+W}{n+p+2}
  \right) \le 
  \xi(W)
  \le 
  \frac{1}{ng_1(W)}, 
\end{equation}
is a necessary condition on 
$\hat{\bm{M}}(\xi,\eta;\bdeltap)$ 
to be admissible in the class (\ref{class:matrix}).

We can easily see that $\hat{\bm{M}}(\xi_0,\eta_0;\bdeltap)$ 
is always nonnegative definite. 
However $l_1(\xi_0;\bdeltap)=0$ when $\xi_0(W) = 1/(ng(W))$
and $l_0(\eta_0;\bdeltap)=0$ when 
$\eta_0(W) = (1/n+g_3(W))/g_1(W)$, 
i.e. 
$\hat{\bm{M}}(\xi_0,\eta_0;\bdeltap)$ is not always positive 
definite.
As mentioned in the previous section, from a practical viewpoint, 
the inverse matrices of estimators of the MSE matrix are required 
to form a confidence set as (\ref{ellipsoid}).
In this sense $\hat{\bm{M}}(\xi_0,\eta_0;\bdeltap)$ may still be 
undesirable. In the next section we provide some positive definite
estimators by considering the improvement of the UMVUE of the risk
reduction matrix $\bm{M}^*(\bdeltap)$.

\subsection{An improved estimation of the risk reduction matrix and 
  positive-definite estimation of the MSE matrix} 
\label{sec:4.2}
\begin{figure}[htbp]
 \centering
 \begin{tabular}{cc}
 \includegraphics[scale=0.35]{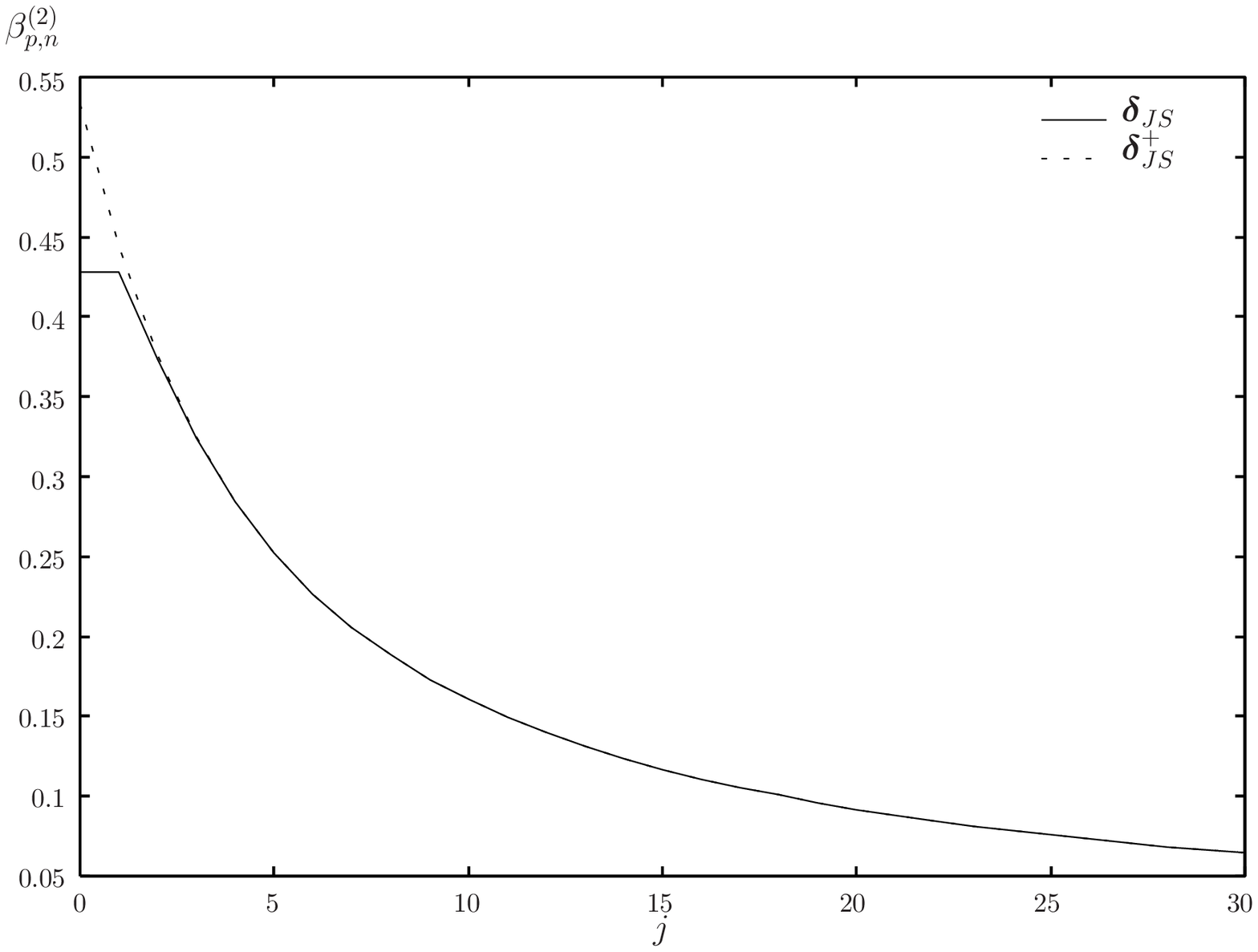}&
 \includegraphics[scale=0.35]{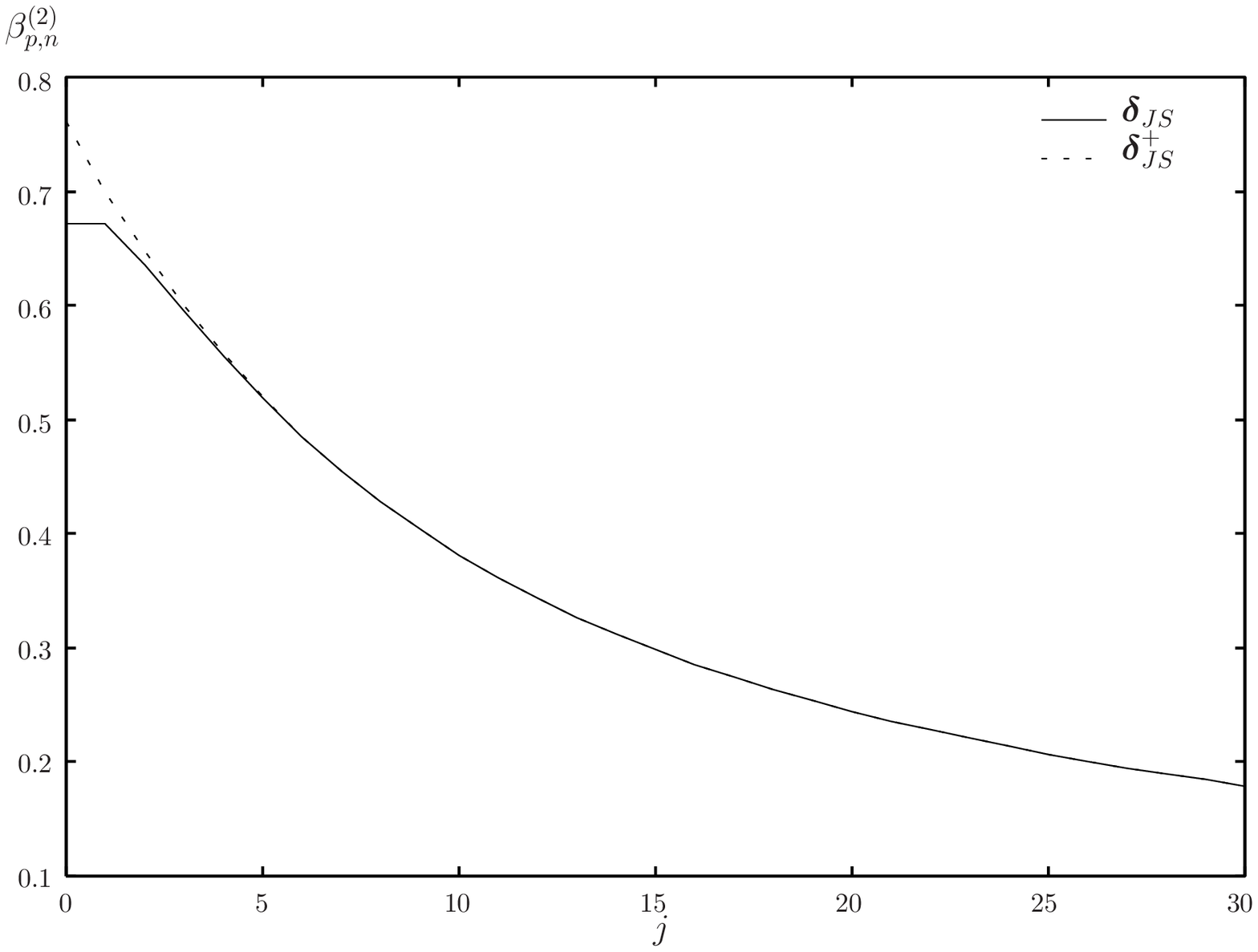}\\
 (i) $(p,n)=(5,5)$ &
 (ii) $(p,n)=(10,10)$ 
 \end{tabular}
 \caption{$\beta^{(2)}_{p,n}(j:\bdeltap)$ of $\bdelta_{JS}$ and
 $\bdelta_{JS}^+$ with $(p,n)=(5,5)$ and $(10,10)$}
\label{figure:beta-2}
\end{figure}

Following the argument in Kubokawa and
Srivastava\cite{Kubokawa-Srivastava},  
we consider the estimation of MSE reduction matrices 
$\bm{M}^*(\bdeltap) = \sigma^2 \bm{I}_p - \bm{M}(\bdeltap)$.
As a criterion, we use the following quadratic loss function, 
\begin{equation}
 \label{loss:reduction-matrix}
  L(\hat{\bm{M}}^*(\bdeltap) ; \bm{M}^*(\bdeltap))
  =
  \mathrm{tr}(\hat{\bm{M}}^*(\bdeltap) - \bm{M}^*(\bdeltap))^2.
\end{equation}
We consider the following class derived from (\ref{class:matrix}), 
\begin{align}
 \label{class:red-matrix}
 \hat{\bm{M}}^*(\xi,\eta ; \bm{\delta}_{\bm{\phi}})
 & = \frac{S}{n}\bm{I}_p - 
 \hat{\bm{M}}(\xi,\eta ; \bm{\delta}_{\bm{\phi}})\notag\\
 & =
 S \cdot \bm{\Gamma}'
 \bigl\{
 g_1(W)\xi(W) \cdot (\bm{I}_p - \bm{E}_{11}) -
 \left(
 g_3(W) - g_1(W)\eta(W)
 \right)
 \cdot \bm{E}_{11}
 \bigr\}
 \bm{\Gamma}\notag\\
 &= S \cdot \bm{\Gamma}' \bm{L}^*(\xi,\eta;\bdeltap) \bm{\Gamma},
\end{align}
where
$$
\bm{L}^*(\xi,\eta;\bdeltap) = 
\mathrm{diag}
(l^*_0(\eta;\bdeltap),l^*_1(\xi;\bdeltap), \ldots, l^*_1(\xi;\bdeltap)), 
$$
$$
l^*_0(\eta;\bdeltap) = g_1(W)\eta(W) - g_3(W)
= \frac{1}{n} - l_0(\eta;\bdeltap),
$$
$$
l^*_1(\xi;\bdeltap) = g_1(W)\xi(W) = 
\frac{1}{n} - l_1(\xi;\bdeltap)
$$
Let $u$, $v$ and $b(W)$ be defined as in the previous section.
Define $\beta^{(2)}_{p,n}(\bdeltap)$ by 
$$
\beta^{(2)}_{p,n}(\bdeltap) :=
\sup_{j \ge 0}\beta^{(2)}_{p,n}(j;\bdeltap), 
$$
where 
$$
\beta^{(2)}_{p,n}(j;\bdeltap) =
\mathrm{E}\left[
\frac{2\phi(u/v)}{u/v} - 
\frac{b(u/v)}{p+2j}
\right].
$$
Figure \ref{figure:beta-2} represents the behavior of 
$\beta^{(2)}_{p,n}(j;\bdelta_{JS})$ and
$\beta^{(2)}_{p,n}(j;\bdelta_{JS}^+)$  
with $(p,n)=(5,5)$, $(10,10)$ computed by Monte Carlo simulation.  
Table \ref{table:beta-2} presents 
$\beta^{(2)}_{p,n}(\bdelta_{JS})$ and
$\beta^{(2)}_{p,n}(\bdelta_{JS}^+)$ when 
$(p,n)$ is $(5,5)$, $(5,10)$, $(10,5)$ and $(10,10)$. 
Kubokawa and Srivastava\cite{Kubokawa-Srivastava} proved that 
$\sup_{j \ge 0}\beta^{(2)}_{p,n}(j;\bdelta_{JS})$ are attained at
$j=0,1$. 
They also showed numerically that 
$\sup_{j \ge 0}\beta^{(2)}_{p,n}(j;\bdelta_{JS}^+)$ are attained at
$j=0$. 
The result of Figure \ref{figure:beta-2} is consistent to their
results. 

Similar to the argument in Theorem \ref{th:matrix-1}, we can obtain the
following result. 

\begin{table}[htbp]
 \centering
 \caption{$\beta^{(2)}_{p,n}(\bdelta_{JS})$ and
 $\beta^{(2)}_{p,n}(\bdelta_{JS}^+)$ for some $p$ and $n$} 
 \label{table:beta-2}
 \begin{tabular}{ccccc}\hline
  $(p,n)$ & $(5,5)$ & $(10,5)$ & $(5,10)$ & $(10,10)$\\ \hline
  $\beta^{(2)}_{p,n}(\bdelta_{JS})$  & 
 0.4260 & 0.5704 & 0.5008 & 0.6718 \\
  $\beta^{(2)}_{p,n}(\bdelta_{JS}^+)$ &
 0.5332 & 0.6611 & 0.5173 & 0.7616\\ \hline
 \end{tabular}
\end{table}

\begin{theorem}
 \label{th:matrix-2}
 Let the solutions of the equations
 $$
 \frac{1+W}{n+p+2} \cdot \frac{\beta_{p,n}^{(2)}(\bdeltap)}{g_1(W)}
 = 1, \quad
 \frac{1}{g_1(W)} 
 \left(
 g_3(W) + \frac{1+W}{n+p+2} \cdot \beta^{{(2)}}_{p,n}(\bdeltap)
 \right)
 =1
 $$
 be denoted by $W_{p,n}^{\xi}(\bdeltap)$ and 
$W_{p,n}^{\eta}(\bdeltap)$,
 respectively. 
 Define $\xi_1(W)$ and $\eta_1(W)$ by 
 $$
 \xi_1(W) = 
 \min\left(1, \;
 \frac{1+W^{\xi}_{n,p}(\bdeltap)}{n+p+2} 
 \cdot \frac{\beta_{p,n}^{(2)}(\bdeltap)}{g_1(W)}
 \right), 
 $$
 $$
 \eta_1(W) = 
 \min \left[1, \;
 \frac{1}{g_1(W)} 
 \left(
 g_3(W) + \frac{1+W^{\eta}_{n,p}(\bdeltap)}{n+p+2} 
 \cdot \beta^{{(2)}}_{p,n}(\bdeltap)
 \right)
 \right].
 $$
 Then $\hat{\bm{M}}^*(\xi_1,\eta_1; \bdeltap)$ improves on 
 $\hat{\bm{M}}_0^*(\bdeltap)$ under the loss
 (\ref{loss:reduction-matrix}). 
\end{theorem}

The proof of this theorem is given in the Appendix.
By using the similar procedure to Kubokawa and
Srivastava\cite{Kubokawa-Srivastava}, we can also obtain the following
theorem.

\begin{theorem}
 \label{th:matrix-3}
 Suppose that 
 \begin{enumerate}
  \item[(i)] $g_1(W)$ is nonincreasing and $g_3(W)/g_1(W)$ is
	     nondecreasing ;
  \item[(ii)] $\xi(\cdot)$ and $\eta(\cdot)$ are nondecreasing and satisfy
	     $$
	     \lim_{W \to \infty} \xi(W)  = 1, \quad    
	     \lim_{W \to \infty} \eta(W)  = 1;
	     $$
  \item[(iii)] $\xi(W) \ge \xi_2(W)$ and $\eta(W) \ge \eta_2(W)$ where
	     $$
	     \xi_2(W) = 
	     \min\left(
	     1, \; \frac{1}{g_1(W)} \cdot
	     \frac{\beta_{p,n}^{(2)}(\bdeltap)}{n+2} 
	     \right), 
	     $$
	     $$
	     \eta_2(W) =
	     \min \left(
	     1,\; 
	     \frac{1}{g_1(W)}
	     \left(
	     g_3(W) + \frac{\beta_{p,n}^{(2)}(\bdeltap)}{n+2}
	     \right)
	     \right).
	     $$
 \end{enumerate}
  Then $\hat{\bm{M}}^*(\xi,\eta; \bdeltap)$ improves on 
  $\hat{\bm{M}}_0^*(\bdeltap)$ under the loss
  (\ref{loss:reduction-matrix}). 
\end{theorem}

The proof of this theorem is also given in the Appendix.
We can easily see from (\ref{g:JS}) and (\ref{g:JS^+})
that the condition (i) is satisfied for 
$\bdelta_{JS}$ and $\bdelta_{JS}^+$.
We note that if $g_3(W) - g_1(W) \ge 0$ for all $W \ge 0$, 
$\eta_1(W)=1$ and $\eta_2(W)=1$ with probability one.

The results of Theorem \ref{th:matrix-2} and Theorem \ref{th:matrix-3} 
lead us to propose the estimators of the MSE matrix $\bm{M}(\bdeltap)$
as 
$$
\hat{\bm{M}}(\xi_i,\eta_i ;\bm{\delta}_{\bm{\phi}})
=
\frac{S}{n}\bm{I}_p - \hat{\bm{M}}^*(\xi_i,\eta_i ;\bm{\delta}_{\bm{\phi}}),
\quad
i=1,2.
$$
If $(p,n)$ satisfies
\begin{align*}
  \gamma_{p,n}^{\xi}(\bdeltap) &:= 
  \frac{n(1+W^{\xi}_{p,n}(\bdeltap))\beta_{p,n}^{(2)}(\bdeltap)}
  {n+p+2}
  \le 1, \\
  \gamma_{p,n}^{\eta}(\bdeltap) &:= 
  \frac{n(1+W^{\eta}_{p,n}(\bdeltap))\beta_{p,n}^{(2)}(\bdeltap)}
  {n+p+2} 
  \le 1, 
\end{align*}
$\hat{\bm{M}}(\xi_1,\eta_1 ;\bm{\delta}_{\bm{\phi}})$ is always positive 
definite and if $(p,n)$ satisfies
$$
\beta^{(2)}_{p,n}(\bdeltap)/(n+2) \le 1/n, 
$$
$\hat{\bm{M}}(\xi_2,\eta_2 ;\bm{\delta}_{\bm{\phi}})$ is always
positive.
Table \ref{table:gamma-xi} present
$\gamma^{\xi}_{p,n}(\bdelta_{JS})$,
$\gamma^{\xi}_{p,n}(\bdelta_{JS}^+)$ and 
$\gamma^{\eta}_{p,n}(\bdelta_{JS}^+)$ when 
$(p,n)$ is $(5,5)$, $(5,10)$, $(10,5)$ and $(10,10)$. 
From Table \ref{table:beta-2} and \ref{table:gamma-xi}, 
we can see that the above conditions are satisfied for $(p,n)$
in the tables.
In Table \ref{table:W-xi-eta} we also presents
$W^{\xi}_{p,n}(\bdelta_{JS})$, 
$W^{\xi}_{p,n}(\bdelta_{JS}^+)$ and 
$W^{\eta}_{p,n}(\bdelta_{JS}^+)$ for the same $(p,n)$ as in
Table \ref{table:gamma-xi}.

\begin{table}[b]
 \centering
 \caption{$\gamma^{\xi}_{p,n}(\bdelta_{JS})$, 
 $\gamma^{\xi}_{p,n}(\bdelta_{JS}^+)$ and 
 $\gamma^{\eta}_{p,n}(\bdelta_{JS}^+)$}
 \label{table:gamma-xi}
 \begin{tabular}{ccccc}\hline
  $(p,n)$ & $(5,5)$ & $(10,5)$ & $(5,10)$ & $(10,10)$\\ \hline
  $\gamma^{\xi}_{p,n}(\bdelta_{JS})$  & 
  0.4312 & 0.6143 & 0.5273 & 0.7548\\
  $\gamma^{\xi}_{p,n}(\bdelta_{JS}^+)$ &
  0.4963 & 0.6690 & 0.6102 & 0.8170 \\
  $\gamma^{\eta}_{p,n}(\bdelta_{JS}^+)$ &
  0.2708 & 0.2567 & 0.4130 & 0.4014\\ \hline
 \end{tabular}
 \centering
 \caption{$W^{\xi}_{p,n}(\bdelta_{JS})$, 
 $W^{\xi}_{p,n}(\bdelta_{JS}^+)$ and 
 $W^{\eta}_{p,n}(\bdelta_{JS}^+)$}
 \label{table:W-xi-eta}
 \begin{tabular}{ccccc}\hline
  $(p,n)$ & $(5,5)$ & $(10,5)$ & $(5,10)$ & $(10,10)$\\ \hline
  $W^{\xi}_{p,n}(\bdelta_{JS})$ &
  1.4198 & 2.6577 & 0.7901 & 1.472\\
  $W^{\xi}_{p,n}(\bdelta_{JS}^+)$ & 
  1.2336 & 2.4405 & 0.6829 & 1.360\\
  $W^{\eta}_{p,n}(\bdelta_{JS}^+)$ &
  0.2185 & 0.3202 & 0.1391 & 0.1596\\ \hline
 \end{tabular}
\end{table}

Recall the fact that (\ref{ineq:xi^*}) is a necessary condition on
$\hat{\bm{M}}(\xi,\eta;\bdeltap)$ to be admissible in the class
(\ref{class:matrix}). 
By using the similar argument in Theorem \ref{th:MSE-3}, we can also
obtain the estimators improving on  
$\hat{\bm{M}}(\xi_i,\eta_i;\bdeltap)$, $i=1,2$
under the loss (\ref{loss:matrix}). 

\begin{theorem}
 \label{th:matrix-5}
 Suppose that $\beta^{(1)}_{p,n}(\bdeltap) \ge 0$. 
 Let $\xi_i^{TR}(W)$ for $i=1,2$ be defined by  
 $$
 \xi^{TR}_i(W) = \max
 \left[
 \xi_i(W), \;
 \frac{1}{g_1(W)}
 \left(
 \frac{1}{n} - 
 \frac{1+W}{n+p+1}
 \right)
 \right]
 $$
 Then if 
 $$
 \left\{ W \mid  \xi^{TR}_i(W) = 
 \frac{1}{g_1(W)}
 \left(
 \frac{1}{n} -
 \frac{1+W}{n+p+1}
 \right)
 \right\} \neq \emptyset, 
 $$
 $\hat{\bm{M}}(\xi^{TR}_i,\eta_i;\bdeltap)$ improves on 
 $\hat{\bm{M}}(\xi_i,\eta_i;\bdeltap)$
 under the loss (\ref{loss:matrix}).
\end{theorem}

The proof of this theorem is easy and omitted. 
It is obvious that 
when 
$\hat{\bm{M}}(\xi_i,\eta_i;\bdeltap)$ is positive definite, 
$\hat{\bm{M}}(\xi^{TR}_i,\eta_i;\bdeltap)$ is also 
positive definite.
The dominance relation between 
$\hat{\bm{M}}(\xi_j^{TR},\eta_j ; \bdeltap)$ and 
$\hat{\bm{M}}_0(\bdeltap)$ under the loss (\ref{loss:matrix})
is interesting.
In the same way as the case of the MSE, however, it seems difficult to 
clarify it at this point. 

\section{Monte Carlo studies}
\subsection{Risk performance of proposed estimators}
\begin{figure}[htbp]
 \centering
 \begin{tabular}{cc}
 \includegraphics[scale=0.3]{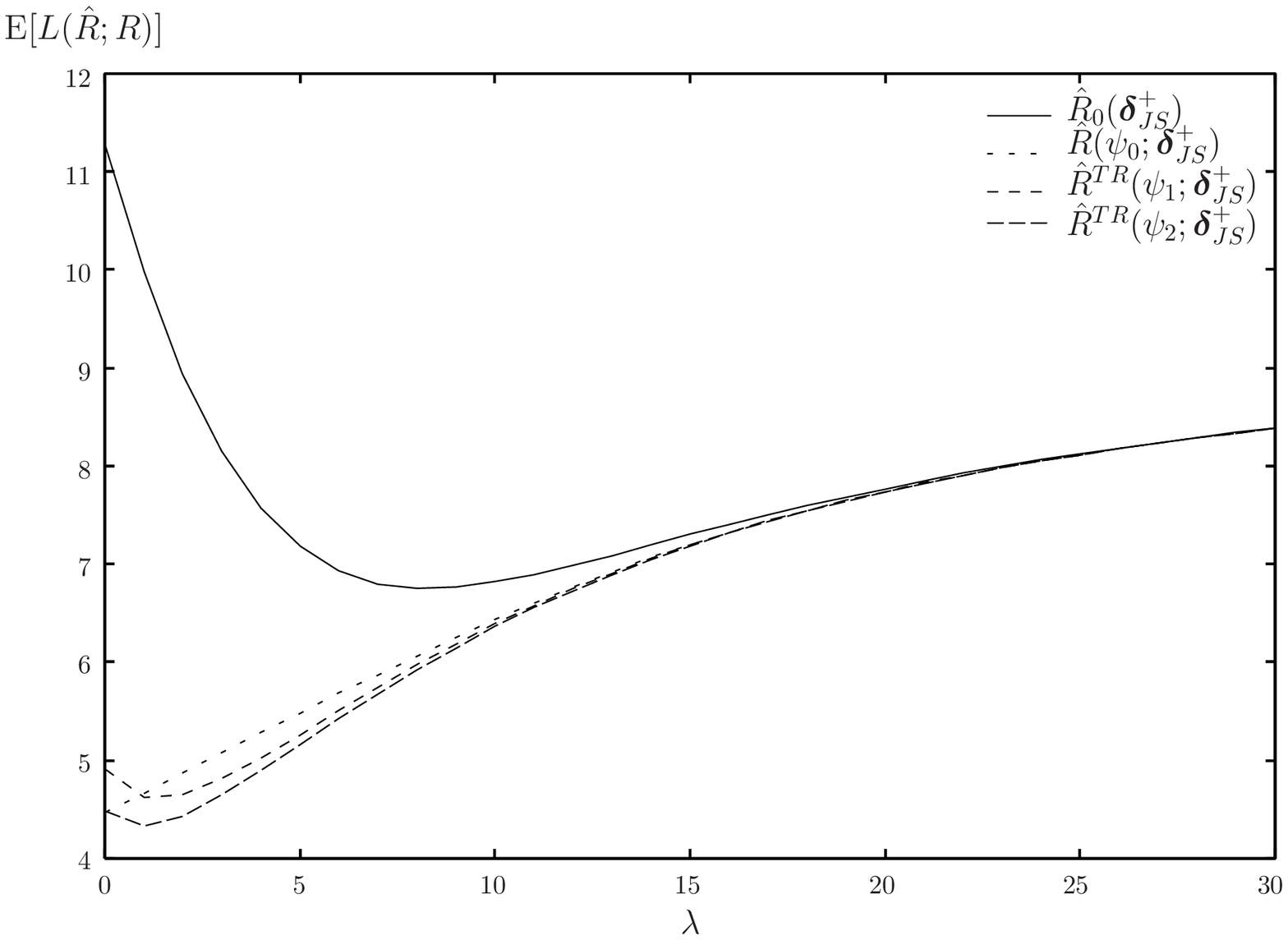} &
 \includegraphics[scale=0.3]{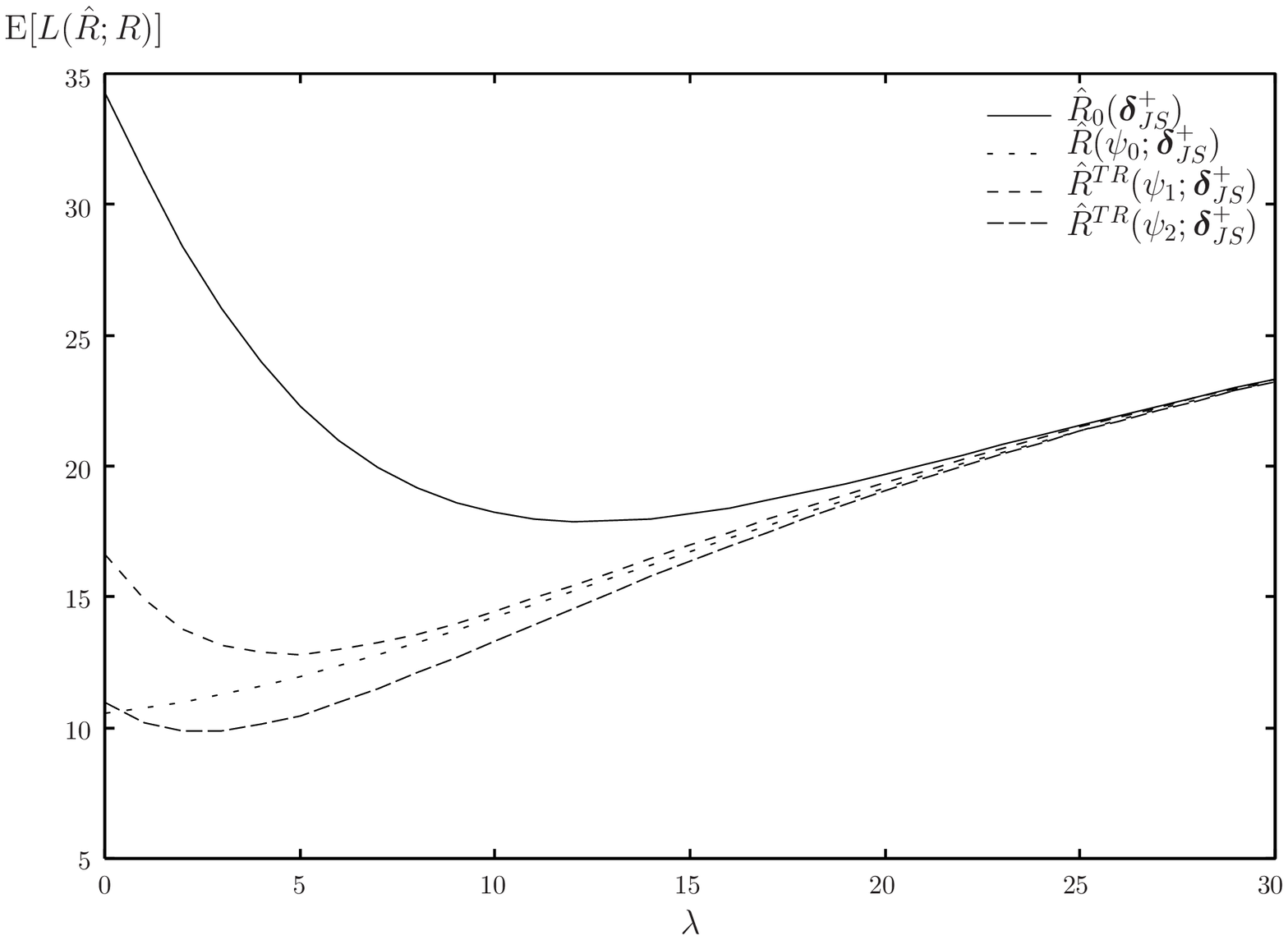}\\
 $(p,n)=(5,5)$ & $(p,n)=(10,5)$\\
 \includegraphics[scale=0.3]{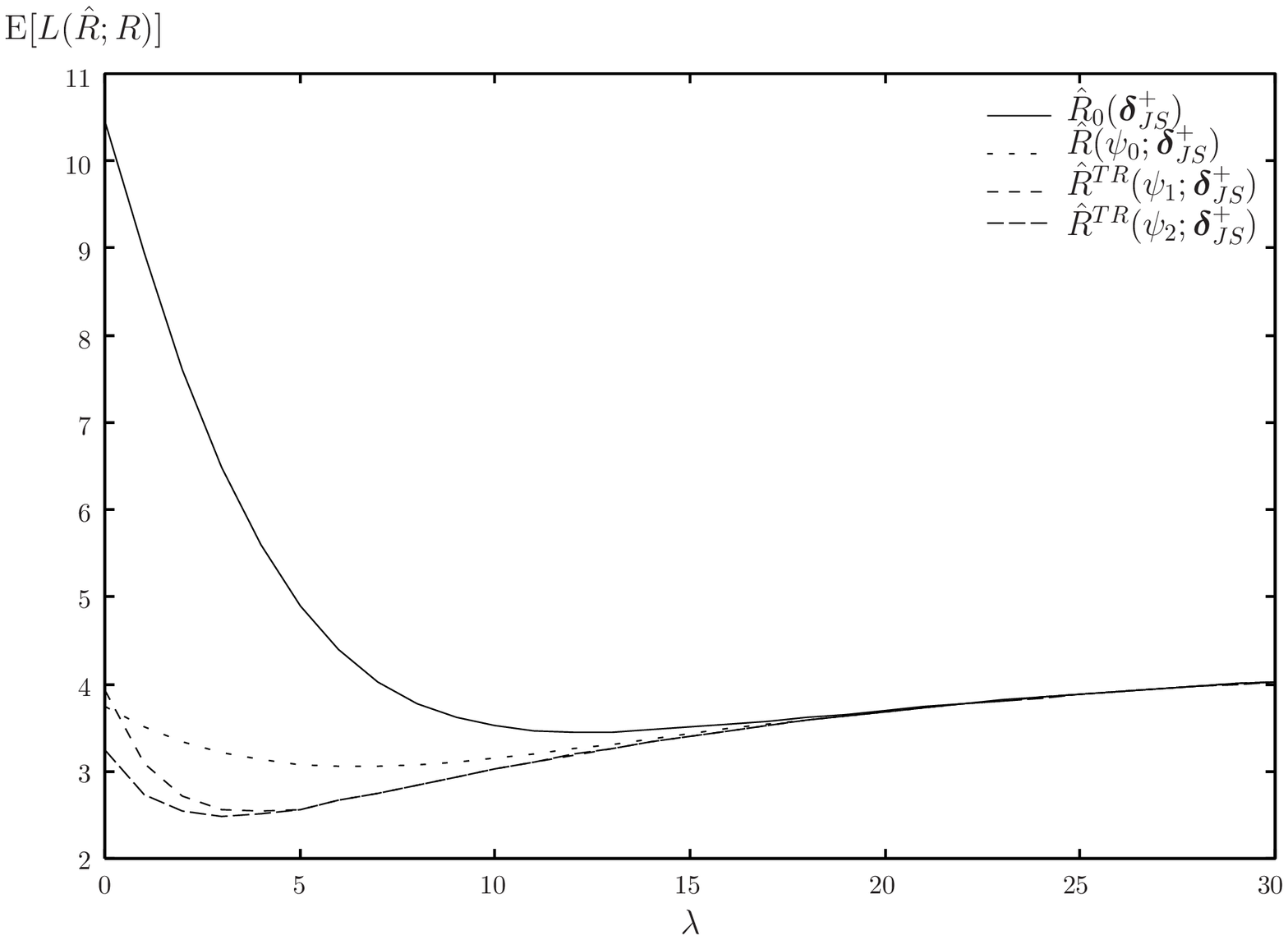} &
 \includegraphics[scale=0.3]{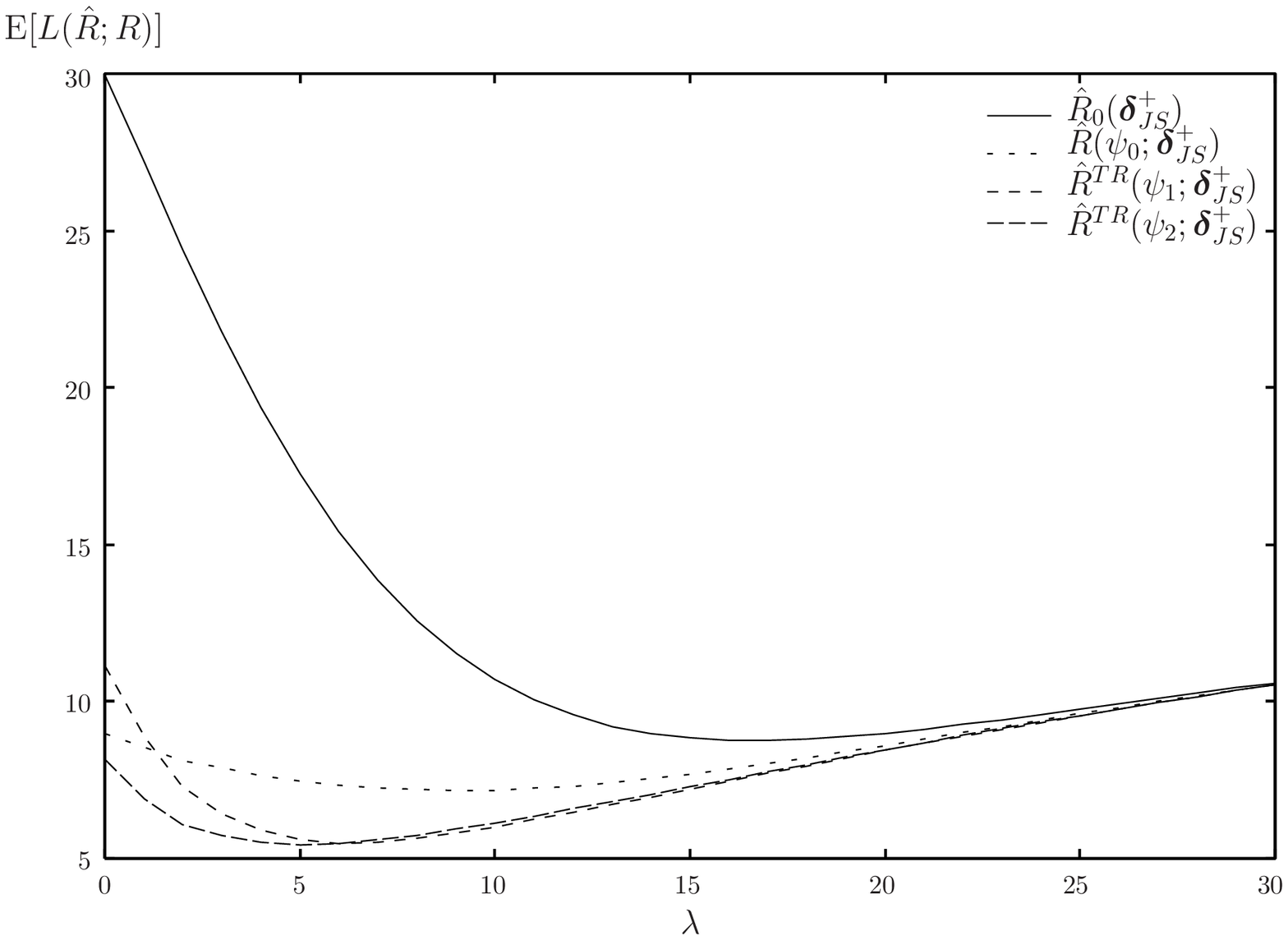}\\
 $(p,n)=(5,10)$ & $(p,n)=(10,10)$
 \end{tabular}
 \caption{Estimated risks of the UMVUE and proposed estimators of 
 $R(\bdelta_{JS}^+)$}
\label{figure:MSE}
\end{figure}

\begin{figure}[htbp]
 \centering
 \begin{tabular}{cc}
 \includegraphics[scale=0.3]{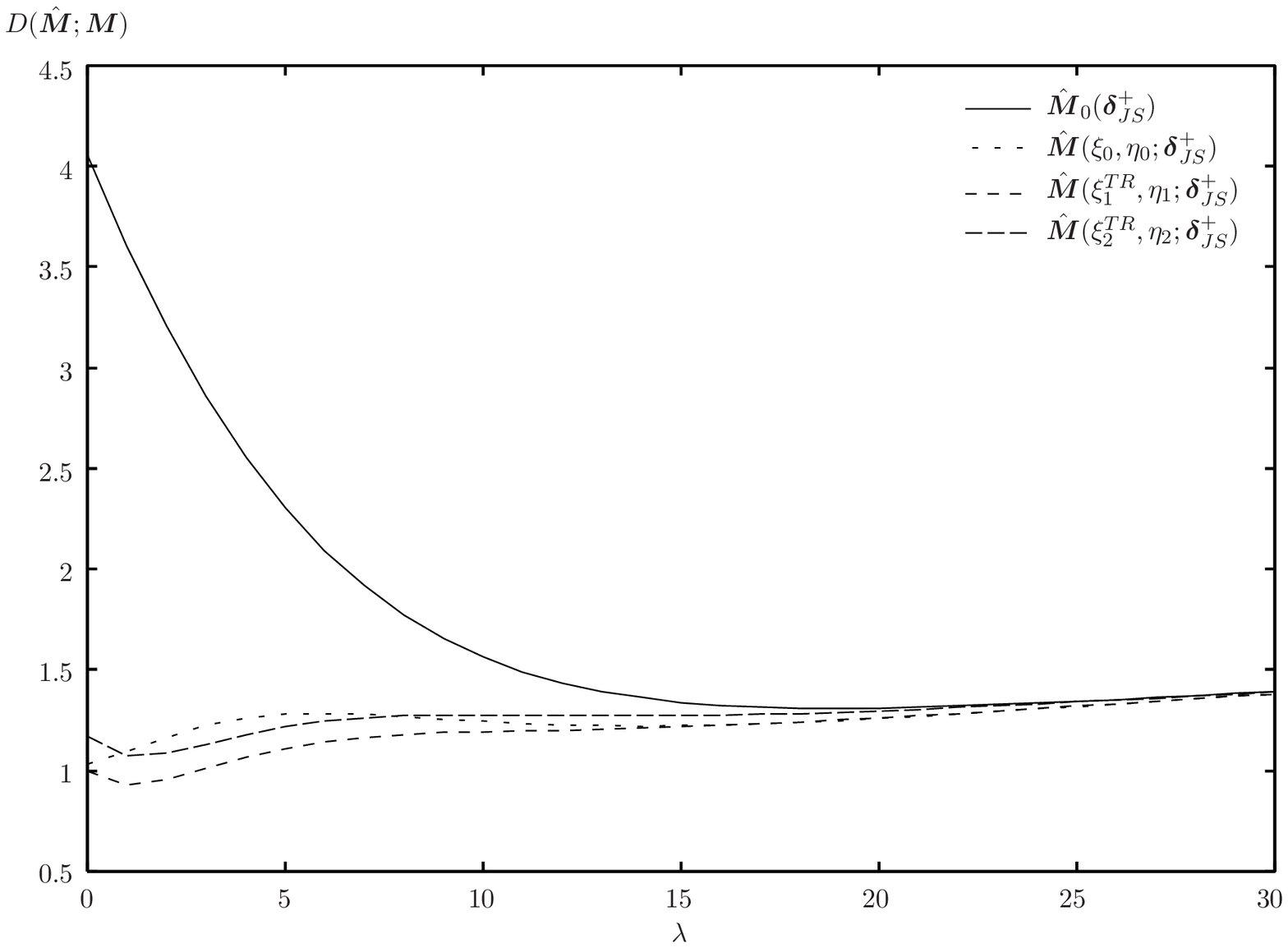} & 
 \includegraphics[scale=0.3]{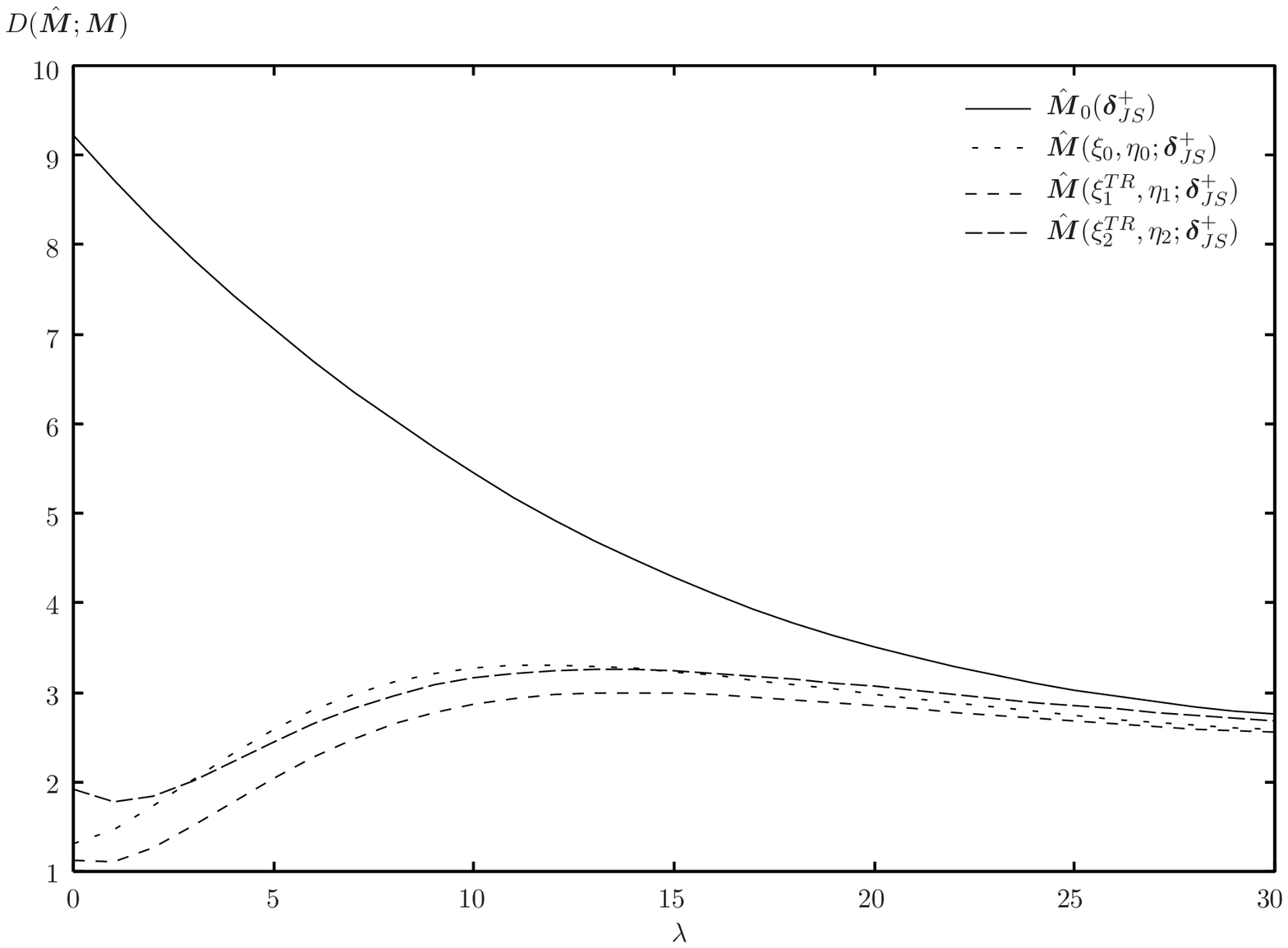}\\
 $(p,n)=(5,5)$ & $(p,n)=(10,5)$\\
 \includegraphics[scale=0.3]{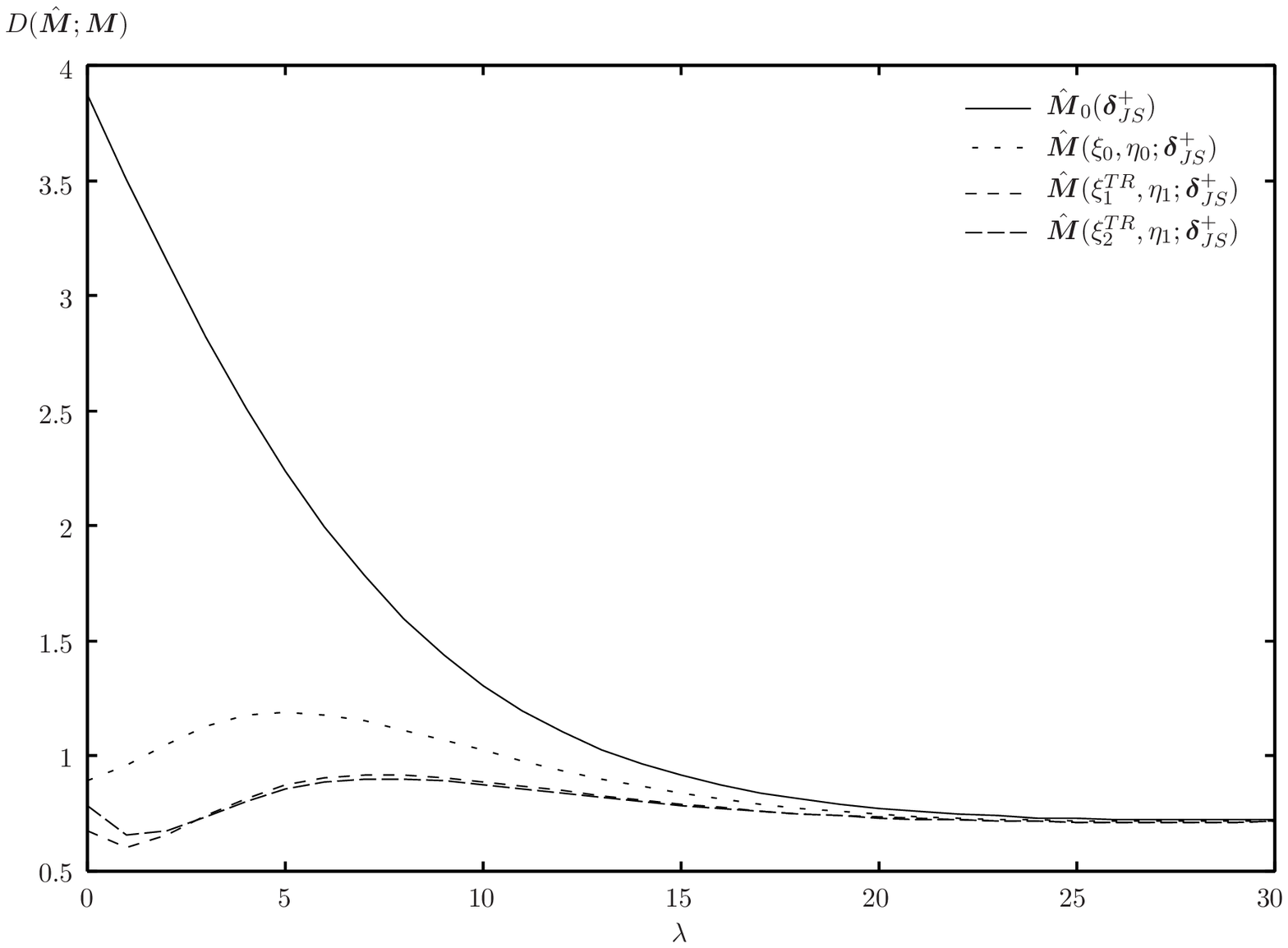} & 
 \includegraphics[scale=0.3]{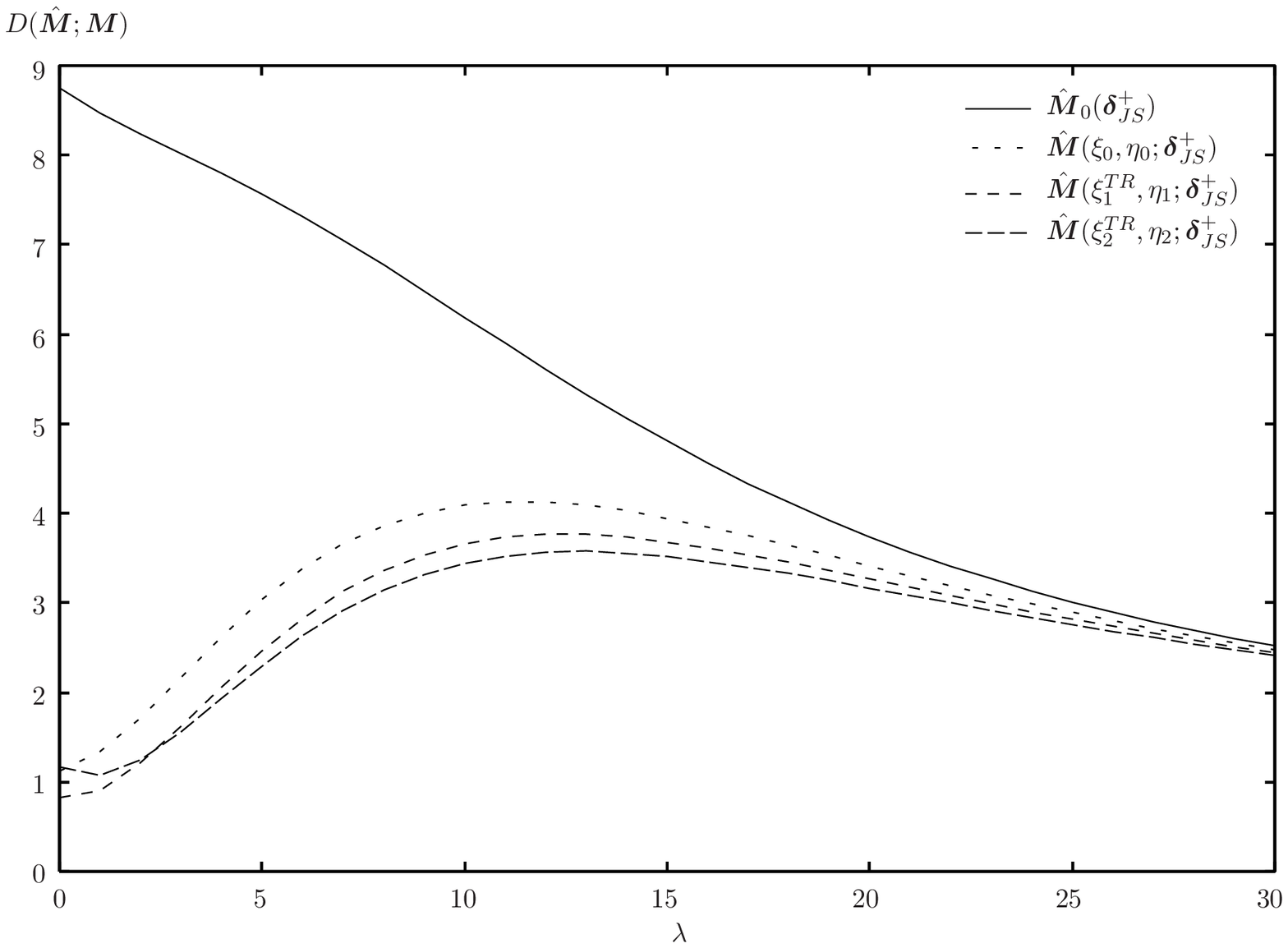}\\
 $(p,n)=(5,10)$ & $(p,n)=(10,10)$
 \end{tabular}
 \caption{Estimated risks for the UMVUE and proposed estimators of 
 $\bm{M}(\bdelta_{JS}^+)$}
\label{figure:matrix-1}
\end{figure}

In this section we study the risk performance of the proposed estimators
of the MSE and the MSE matrix for the positive part Stein estimator
through Monte Carlo studies with 100,000 replications.
Figure \ref{figure:MSE} represents the performance of risks
$\mathrm{E}[L(\hat{R}(\bdelta_{JS}^+) ; R(\bdelta_{JS}^+))]$ of 
$\hat{R}_0(\bdelta_{JS}^+)$, 
$\hat{R}(\psi_0;\bdelta_{JS}^+)$, 
$\hat{R}(\psi_1^{TR};\bdelta_{JS}^+)$ and 
$\hat{R}(\psi_2^{TR};\bdelta_{JS}^+)$. 
Figure \ref{figure:matrix-1} represents the performance of risks 
$\mathrm{E}[L(\hat{\bm{M}}(\bdelta_{JS}^+) ; \bm{M}(\bdelta_{JS}^+))]$
of $\hat{\bm{M}}_0(\bdelta_{JS}^+)$, 
$\hat{\bm{M}}(\psi_0;\bdelta_{JS}^+)$, 
$\hat{\bm{M}}(\psi_1^{TR};\bdelta_{JS}^+)$ and 
$\hat{\bm{M}}(\psi_2^{TR};\bdelta_{JS}^+)$. 
We set $\bm{\theta}$, $\sigma^2$, $\lambda$, $p$ and $n$ as 
\begin{itemize}
 \item $\bm{\theta}=(\sqrt{\lambda/p},\ldots,\sqrt{\lambda/p})'$ ;
 \item $\sigma^2=1$ ;
 \item $0 \le \lambda \le 30$ ;
 \item $(p,n) = (5,5)$, $(10,5)$, $(5,10)$ and $(10,10)$ ;
\end{itemize}

The summary of experiments is as follows.

\begin{itemize}
 \item We can see from the figures that every proposed estimator
       uniformly improves on the UMVUE.
 \item In Section \ref{sec:3}, we have shown that
       $\hat{R}(\psi_0;\bdelta_{JS}^+)$ improves on 
       $\hat{R}_0(\bdelta_{JS}^+)$ under the loss (\ref{loss:MSE}).
       However we could not give theoretical proofs that
       $\hat{R}(\psi_1^{TR};\bdelta_{JS}^+)$ and  
       $\hat{R}(\psi_2^{TR};\bdelta_{JS}^+)$ improve on 
       $\hat{R}_0(\bdelta_{JS}^+)$ under the loss (\ref{loss:MSE}) as
       estimators of the MSE. 
       However we can see from Figure \ref{figure:MSE} that 
       $\hat{R}(\psi_1^{TR};\bdelta_{JS}^+)$ and 
       $\hat{R}(\psi_2^{TR};\bdelta_{JS}^+)$ seem to uniformly improve
       on $\hat{R}_0(\bdelta_{JS}^+)$ under the loss (\ref{loss:MSE}). 
 \item Also in the case of estimators of the MSE matrix, 
       we can see from Figure \ref{figure:matrix-1} that
       $\hat{\bm{M}}(\psi_1^{TR};\bdelta_{JS}^+)$ and 
       $\hat{\bm{M}}(\psi_2^{TR};\bdelta_{JS}^+)$
       which are derived from improved estimators of 
       $\bm{M}^*(\bdelta_{JS}^+)$ seem to improve on 
       $\hat{\bm{M}}_0(\bdelta_{JS}^+)$ under the loss
       (\ref{loss:matrix}).  
 \item When $\lambda$ is small, proposed estimators show large rate of 
       risk reduction. 
       Conversely, 
       as $\lambda$ gets larger, the rate of risk reduction of the
       proposed estimators gets smaller.
\end{itemize}

\subsection{Confidence sets centered at shrinkage estimators}
In this section we study the performance of confidence sets
(\ref{ellipsoid}) centered at the positive part Stein estimator.
Let $c = F_{0.95}(p,n)$ be 95\% percentile of $F$ distribution with 
degrees of freedom $p$ and $n$. 
Then the conventional 95\% confidence set of $\bm{\theta}$ is written
by  
$$
C_0 \; : \; 
Q_0 = 
\frac{(\bm{X} - \bm{\theta})'(\bm{X} - \bm{\theta})/p}{S/n} \le c.  
$$
In what follows, we write 
$\hat{\bm{M}}_1 = \hat{\bm{M}}(\xi_1^{TR},\eta_1;\bdelta_{JS}^+)$ and 
$\hat{\bm{M}}_2 = \hat{\bm{M}}(\xi_2^{TR},\eta_2;\bdelta_{JS}^+)$ for
notational simplicity.
Following Carter et al.\cite{Carter-etal} and Wan et al.\cite{Wan-etal},
we first consider the following three alternatives to $C_0$, 
\begin{align*}
 C_1 \; : \; Q_1 &=
 \frac{(\bdelta_{JS}^+ - \bm{\theta})'\hat{\bm{M}}_1^{-1}
 (\bdelta_{JS}^+ - \bm{\theta})}{p}
 \le c,  \\
 C_2 \; : \; Q_2 &=
 \frac{(\bdelta_{JS}^+ - \bm{\theta})'\hat{\bm{M}}_2^{-1}
 (\bdelta_{JS}^+ - \bm{\theta})}{p}
 \le c,  \\
 C_3 \; : \; Q_3 &=
 \frac{(\bdelta_{JS}^+ - \bm{\theta})'
 (\bdelta_{JS}^+ - \bm{\theta})/p}{S/n}
 \le c. 
\end{align*}
Let $Vol(C_i)$ denote the volume of $C_i$, i.e.
$$
Vol(C_0) = Vol(C_3) = 
\frac{(S/n)^{p/2} \cdot (c \pi)^{p/2}}{\Gamma(p/2 + 1)}, 
$$
$$
Vol(C_1) = \frac{\vert \hat{\bm{M}}_1 \vert^{1/2} (c \pi)^{p/2}}
{\Gamma(p/2 + 1)},\quad 
Vol(C_2) = \frac{\vert \hat{\bm{M}}_2 \vert^{1/2} (c \pi)^{p/2}}
{\Gamma(p/2 + 1)}.
$$
Define $V_1$ and $V_2$ by 
$$
V_1 = \frac{\mathrm{E}[Vol(C_1)]}{\mathrm{E}[Vol(C_0)]}
= \frac{\mathrm{E}[Vol(C_1)]}{\mathrm{E}[Vol(C_3)]}, 
\quad
V_2 = \frac{\mathrm{E}[Vol(C_2)]}{\mathrm{E}[Vol(C_0)]}
= \frac{\mathrm{E}[Vol(C_2)]}{\mathrm{E}[Vol(C_3)]}.
$$
\begin{figure}[htbp]
 \centering
 \includegraphics[scale=0.32]{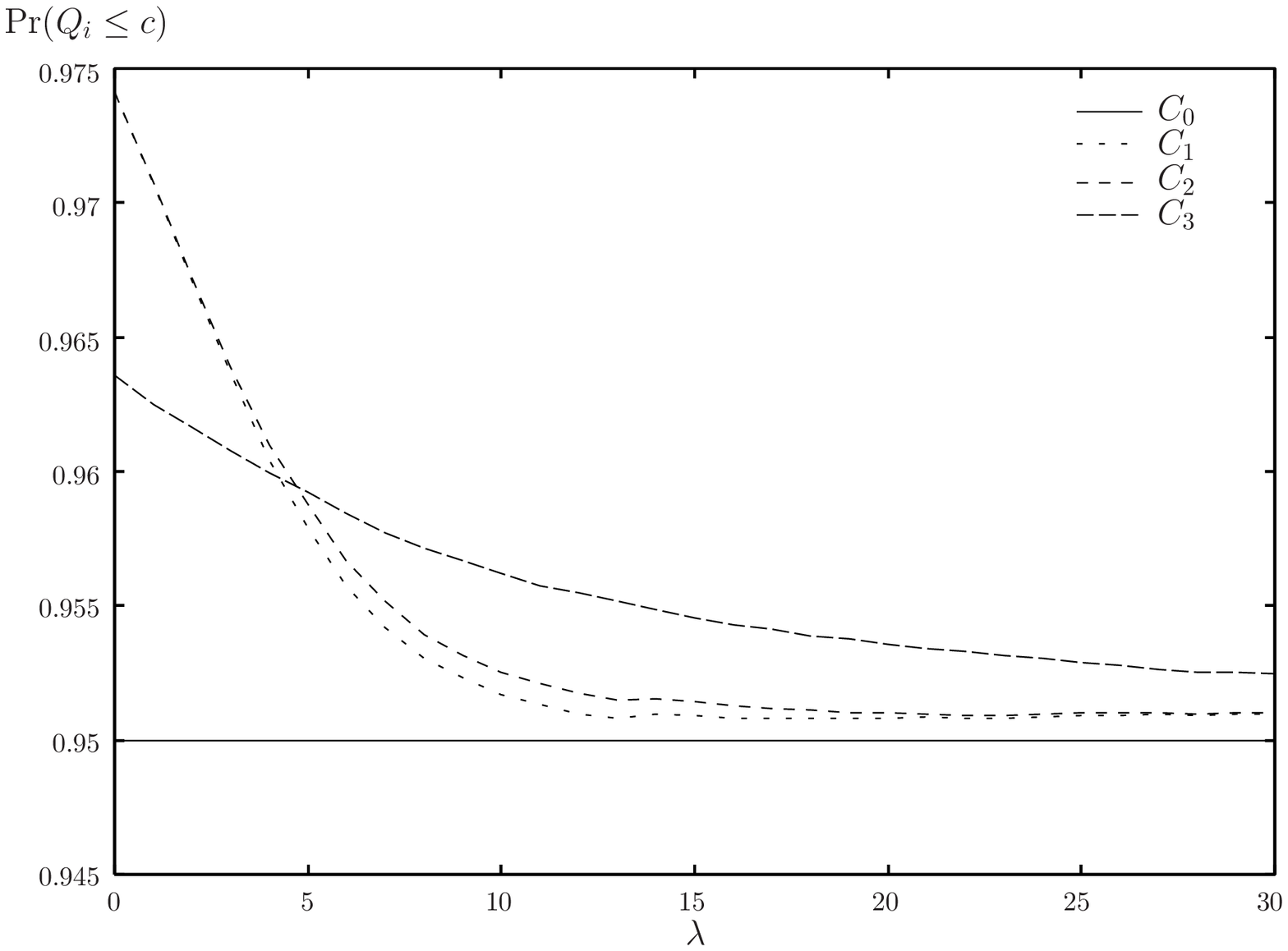}
 \includegraphics[scale=0.32]{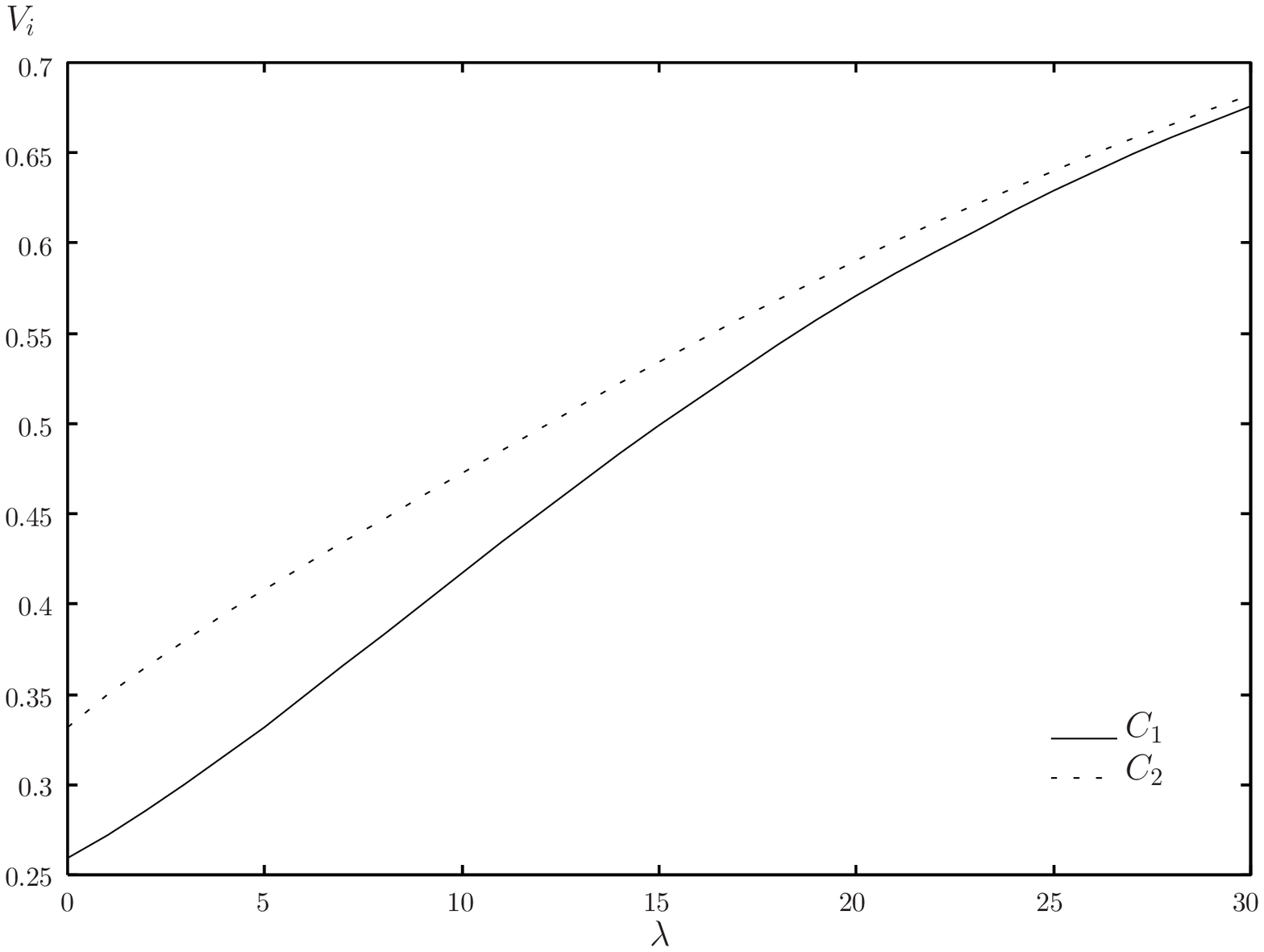}\\
 \hspace*{0.2cm}
 $\mathrm{Pr}(Q_i < c)$ for $(p,n)=(5,5)$ \hspace{1cm}
 $V_1$ and $V_2$ for $(p,n)=(5,5)$\\
 \vspace{0.1cm}
 \includegraphics[scale=0.32]{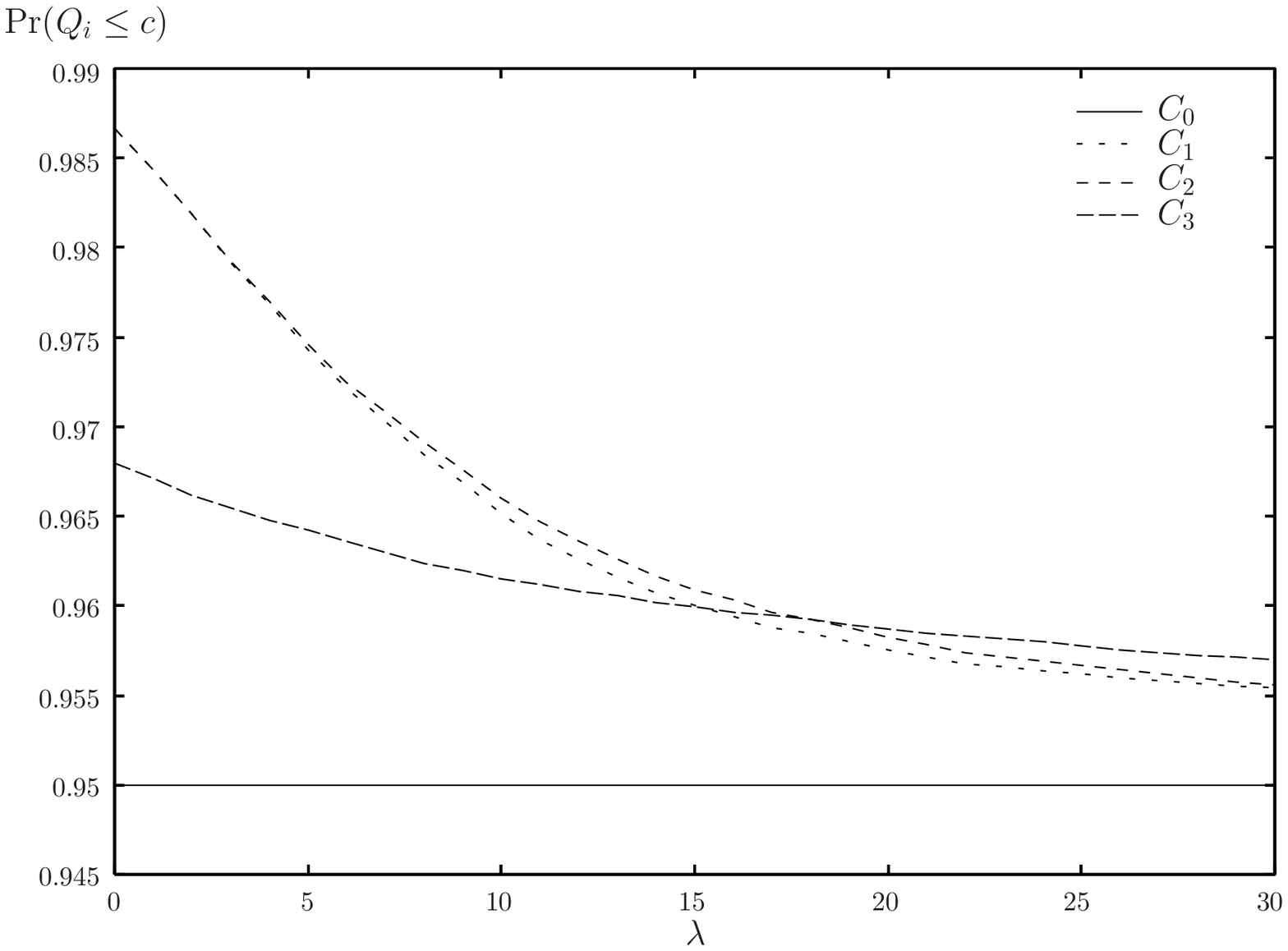}
 \includegraphics[scale=0.32]{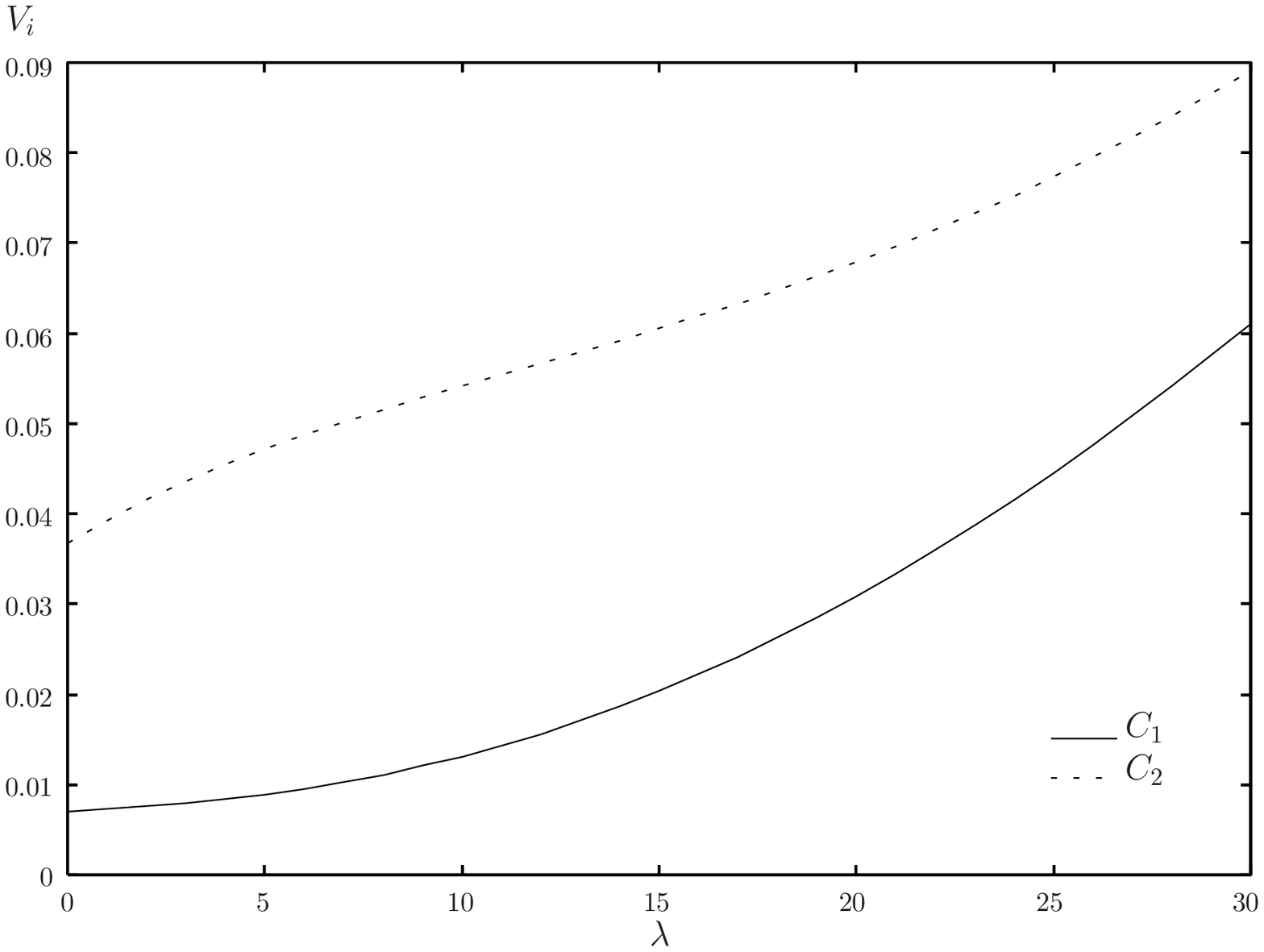}\\
 \hspace*{0.2cm}
 $\mathrm{Pr}(Q_i < c)$ for $(p,n)=(10,5)$ \hspace{1cm}
 $V_1$ and $V_2$ for $(p,n)=(10,5)$\\
 \vspace{0.1cm}
 \includegraphics[scale=0.32]{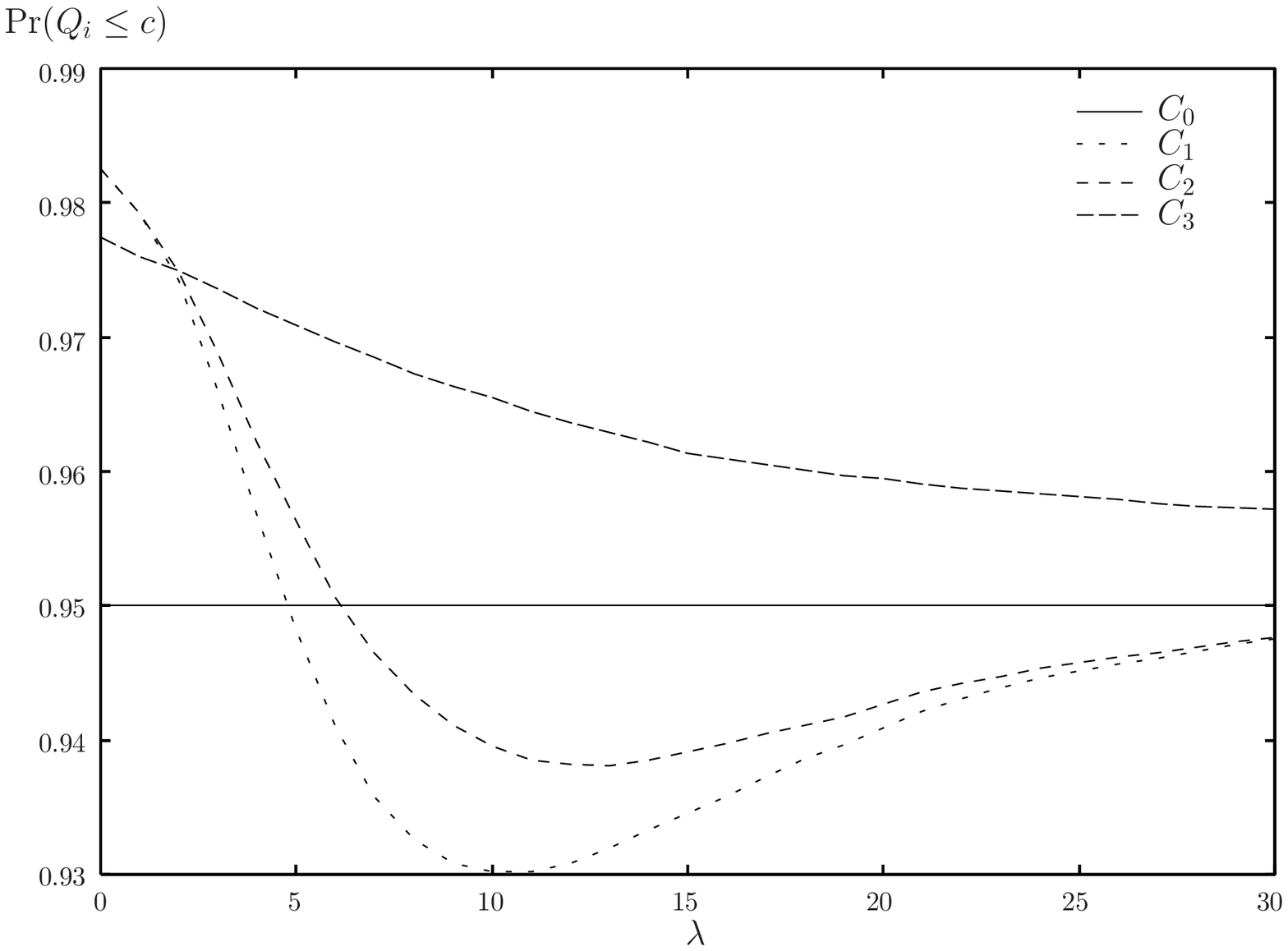}
 \includegraphics[scale=0.32]{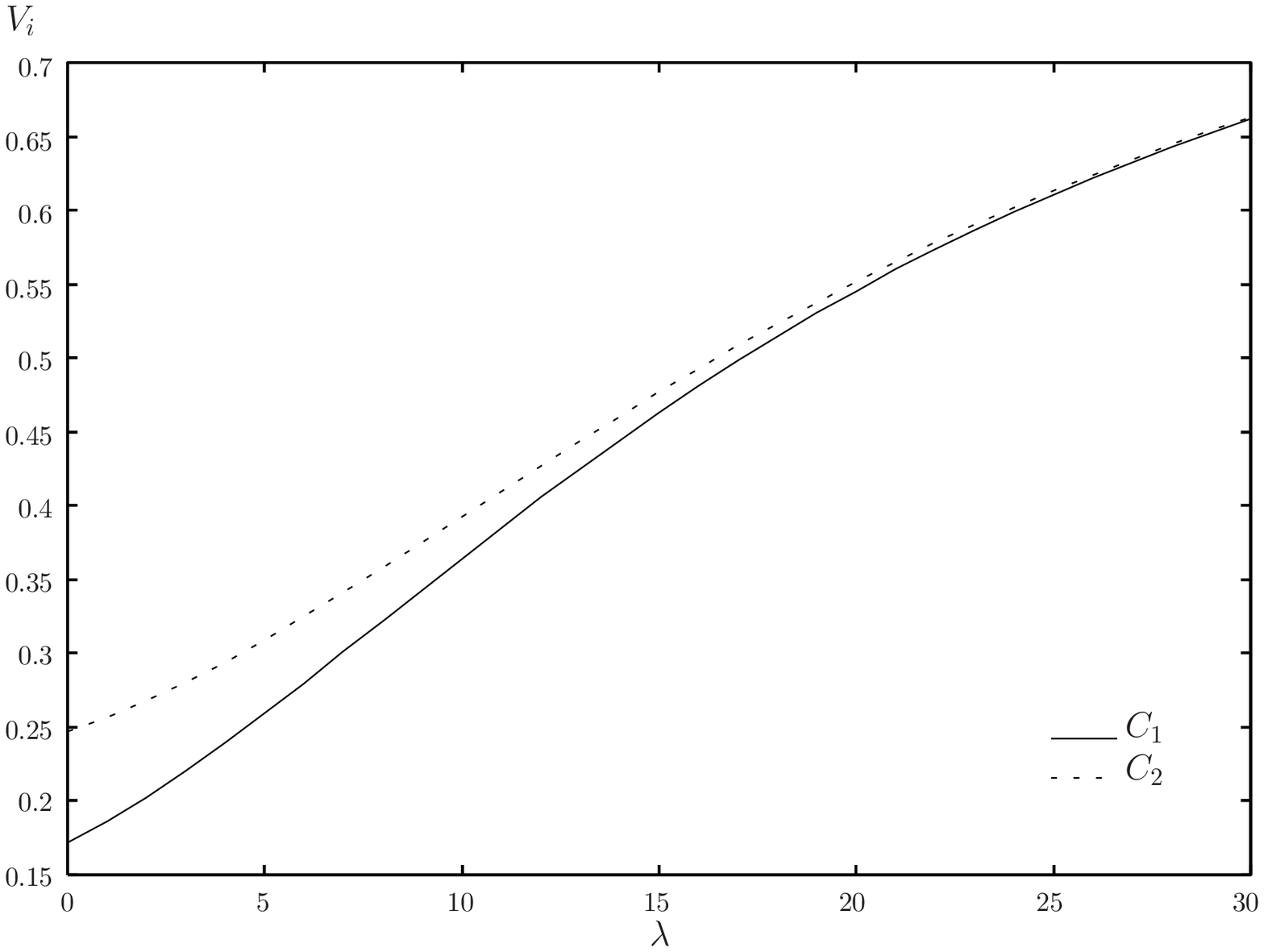}\\
 \hspace*{0.2cm}
 $\mathrm{Pr}(Q_i < c)$ for $(p,n)=(5,10)$ \hspace{1cm}
 $V_1$ and $V_2$ for $(p,n)=(5,10)$\\
 \vspace{0.1cm}
 \includegraphics[scale=0.32]{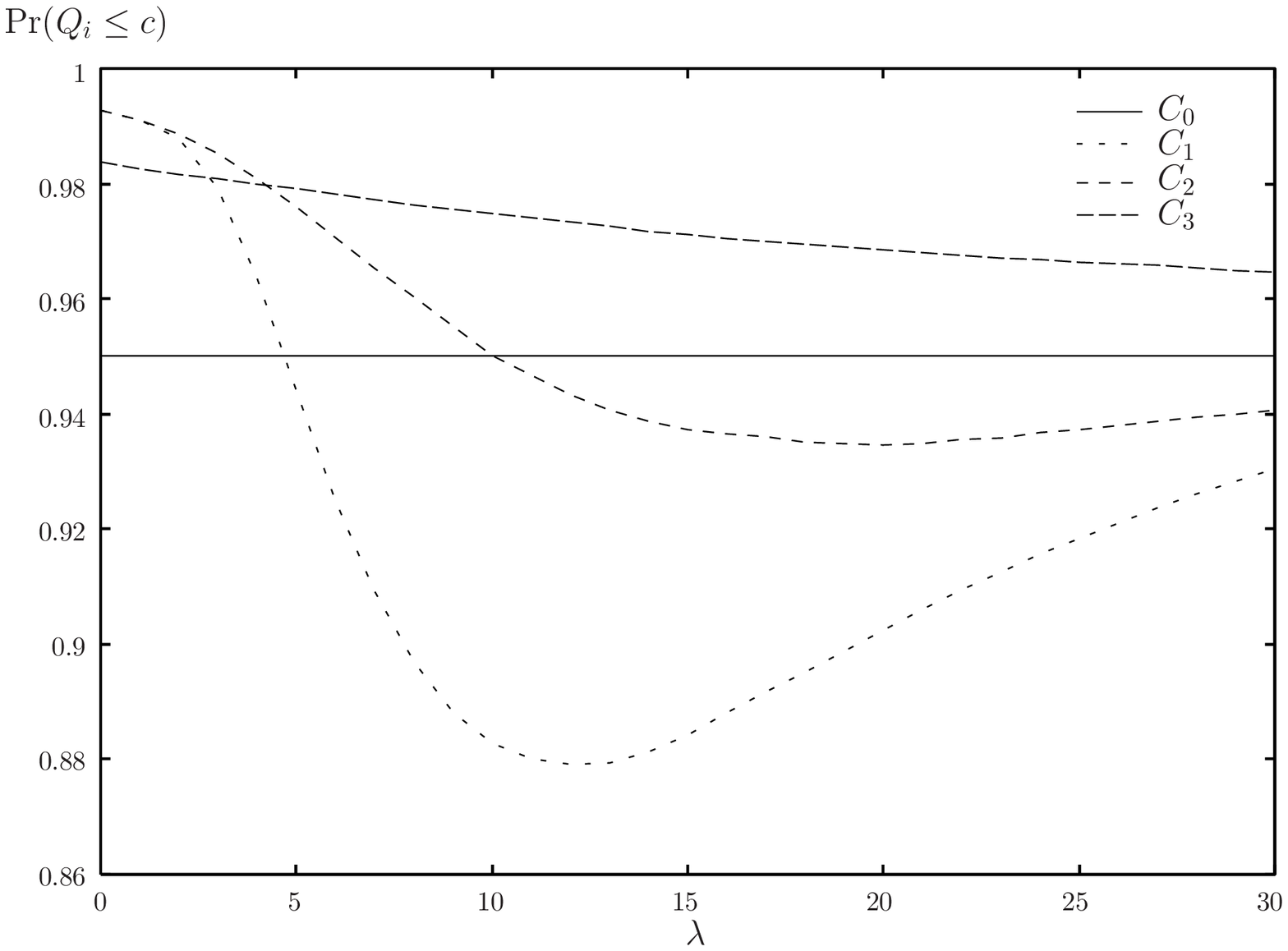}
 \includegraphics[scale=0.32]{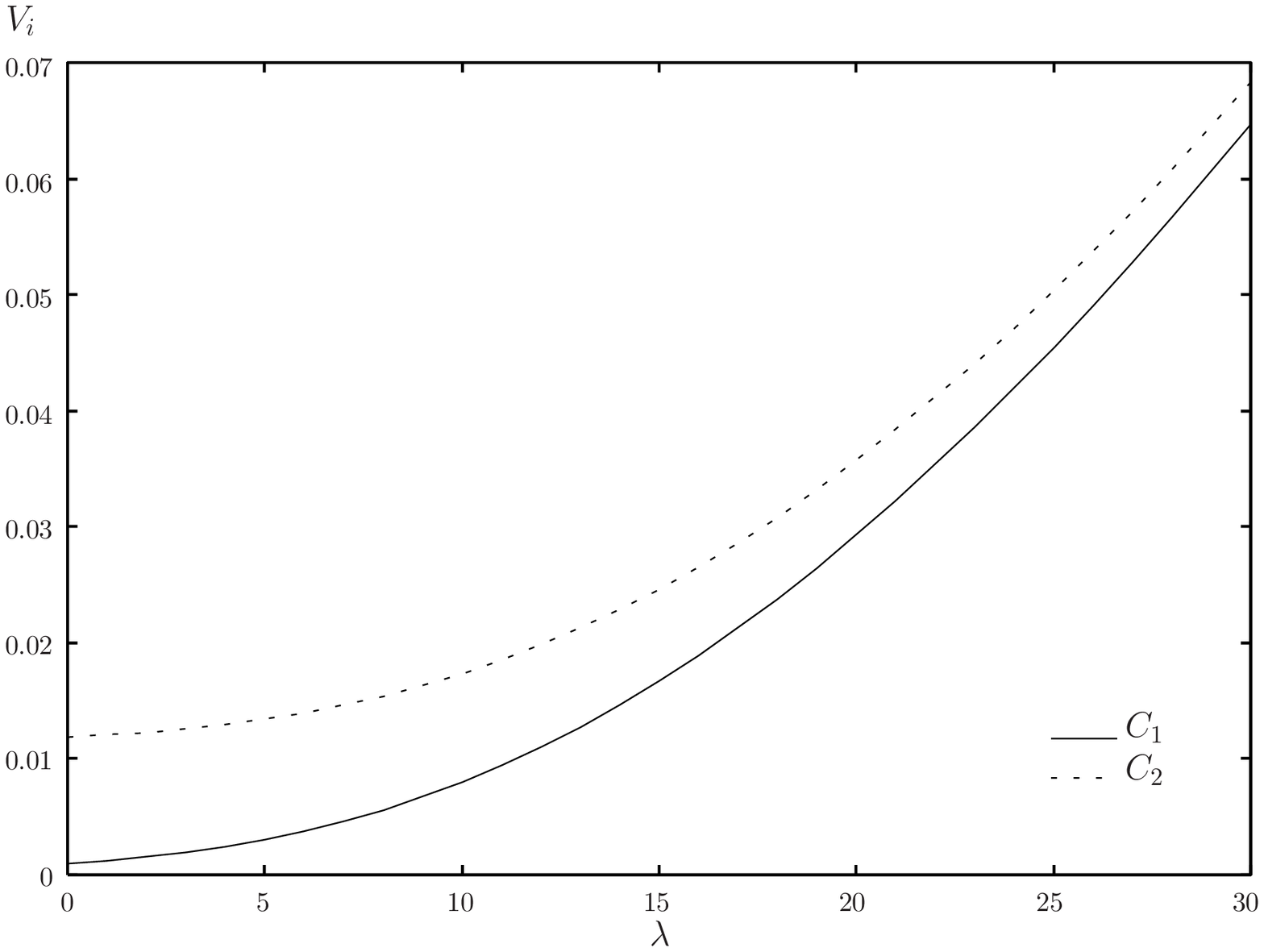}\\
 \hspace*{0.2cm}
 $\mathrm{Pr}(Q_i < c)$ for $(p,n)=(10,10)$ \hspace{1cm}
 $V_1$ and $V_2$ for $(p,n)=(10,10)$\\
 \vspace{0.1cm}
 \caption{Coverage probabilities of $Q_0$ to $Q_3$ and $V_1$ and $V_2$.}
\label{figure:prob_volume-1}
\end{figure}

\begin{figure}[htbp]
 \centering
 \includegraphics[scale=0.3]{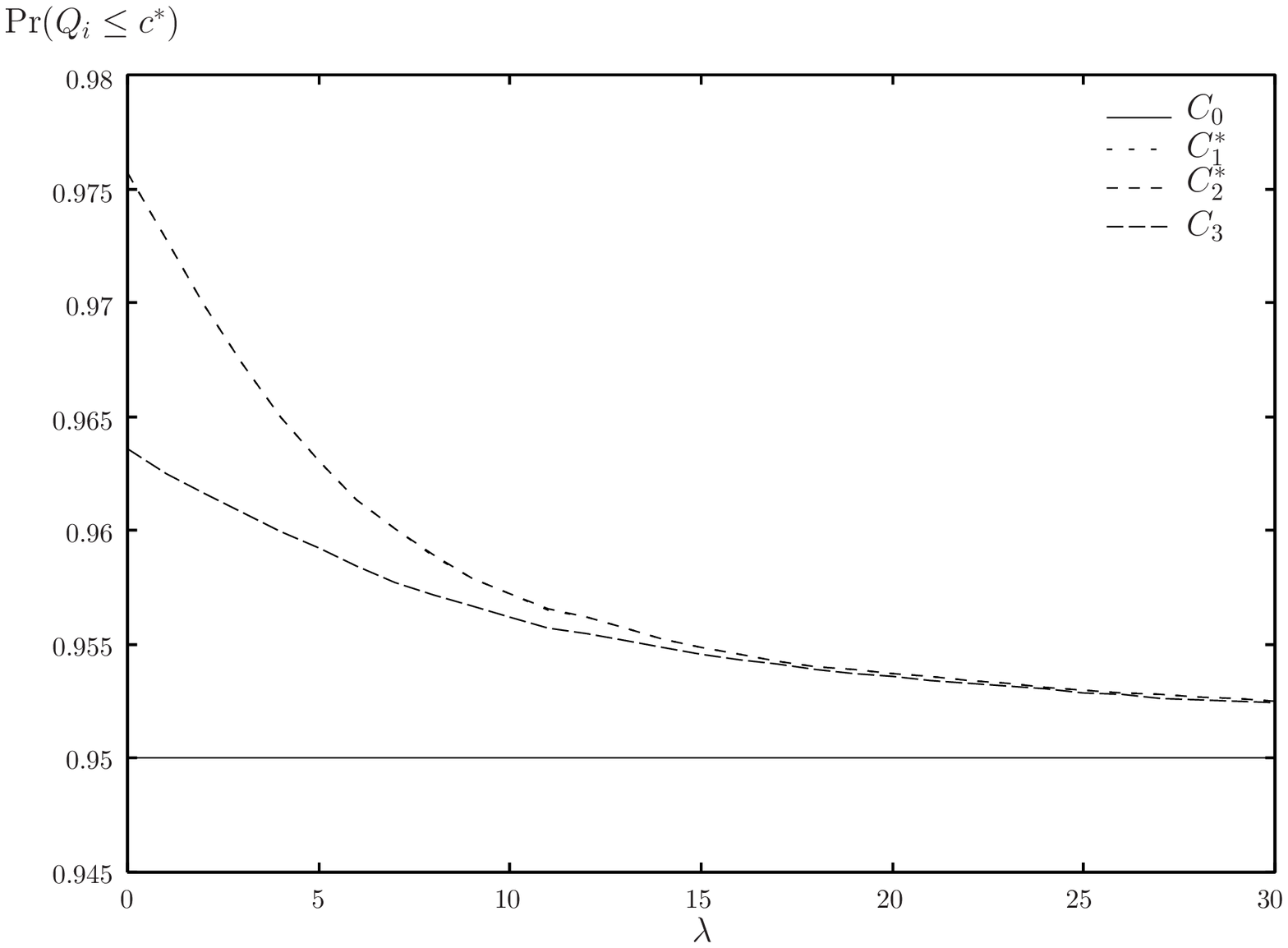}
 \includegraphics[scale=0.3]{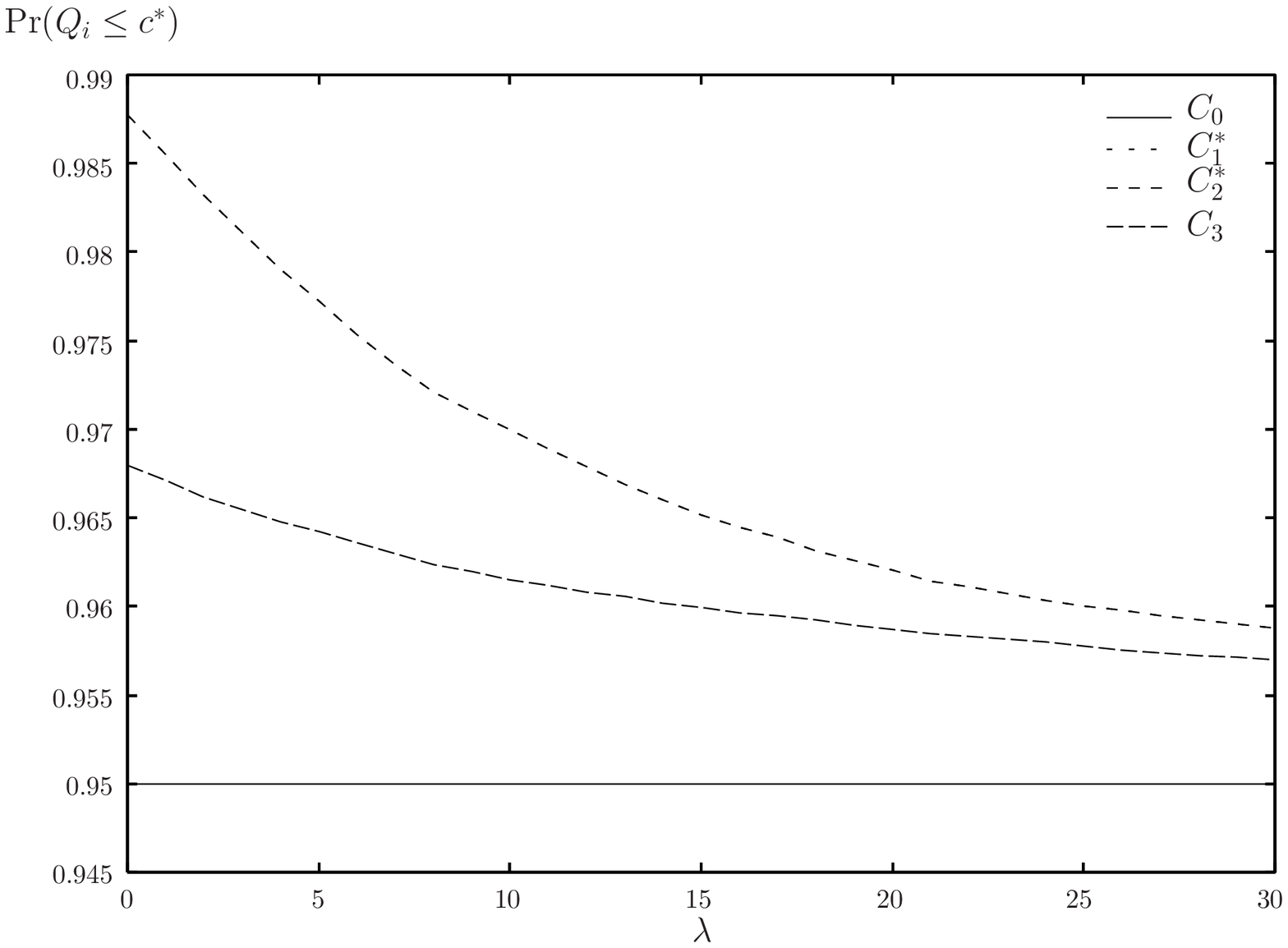}\\
 $(p,n)=(5,5)$\hspace{3.5cm}$(p,n)=(10,5)$\\
 \vspace{0.1cm}
 \includegraphics[scale=0.3]{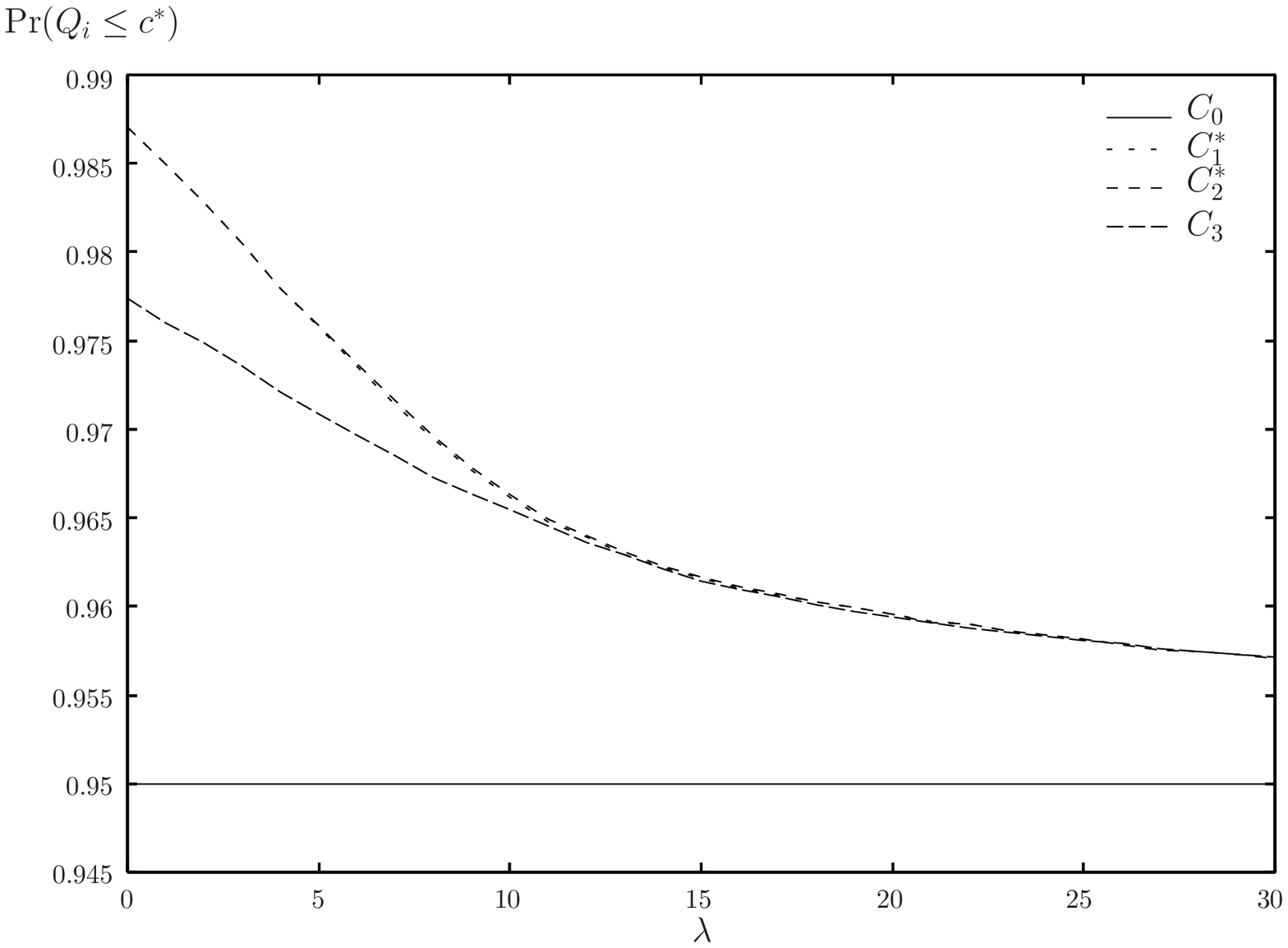}
 \includegraphics[scale=0.3]{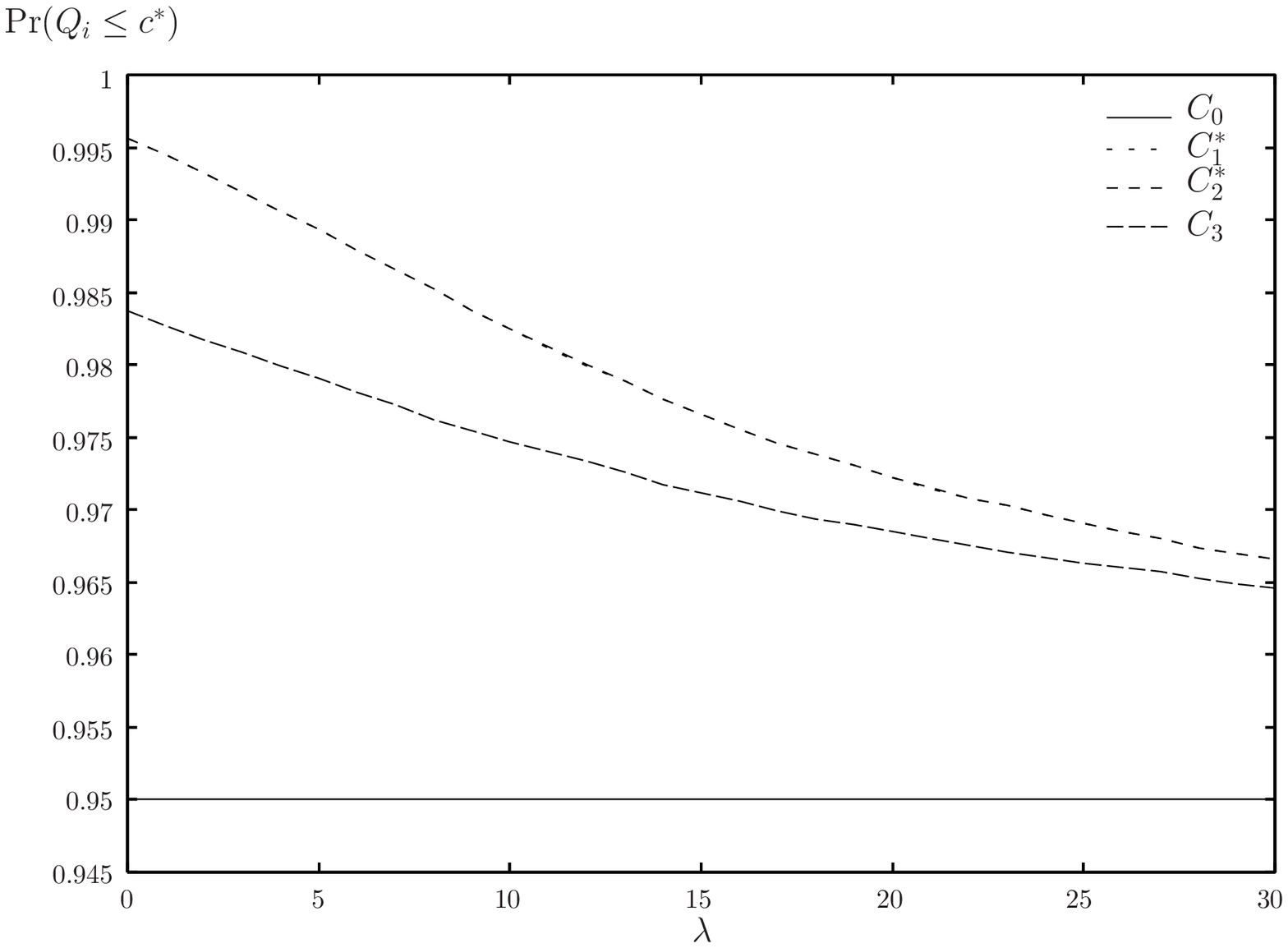}\\
 $(p,n)=(5,10)$\hspace{3.5cm}$(p,n)=(10,10)$\\
 \caption{Coverage probabilities of $C_1^*$ and $C_2^*$}
\label{figure:prob^*-1}
\end{figure}

Figure \ref{figure:prob_volume-1} represent coverage probabilities of
$C_0$ to $C_3$, $V_1$ and $V_2$ computed through Monte Carlo method
with 100,000 replications. 
We set $\bm{\theta}$, $\sigma^2$, $\lambda$, $p$ and $n$ as the
experiments in the previous section.
From the figures, we can see that when $\lambda$ is small, 
the coverage probabilities of $C_1$ to $C_3$ are larger than 
that of $C_0$. 
When $(p,n)$ is $(5,5)$ and $(10,5)$, 
the coverage probabilities of $C_1$ to $C_3$ are
larger than $95\%$ for $0 \le \lambda \le 30$. 
When $(p,n)$ is $(5,10)$ and $(10,10)$, however, 
the coverage probabilities of $C_1$ and $C_2$ drop to less than 
$95\%$ as $\lambda$ gets larger.
The expected volumes of $C_1$ and $C_2$ are far smaller than those of  
$C_0$ and $C_3$.
These results suggest that when $\lambda$ and $n$ is small, 
$C_1$ and $C_2$ seem to have higher coverage probabilities and smaller
volumes than $C_0$. 
On the other hand, while the expected volume of $C_3$ is larger than
that of $C_1$ and $C_2$, the coverage probabilities of $C_3$ is always
larger than $95\%$.   
In this sense, $C_3$ may be the most practical among
$C_0$ to $C_3$.

So we next consider the following confidence sets, 
\begin{align*}
 C_1^* \; : \; Q_1^* &=
 \frac{(\bdelta_{JS}^+ - \bm{\theta})'\hat{\bm{M}}_1^{-1}
 (\bdelta_{JS}^+ - \bm{\theta})}{p}
 \le \frac{(S/n) \cdot c}{\vert \hat{\bm{M}}_1 \vert^{1/p}},  \\
 C_2^* \; : \; Q_2^* &=
 \frac{(\bdelta_{JS}^+ - \bm{\theta})'\hat{\bm{M}}_2^{-1}
 (\bdelta_{JS}^+ - \bm{\theta})}{p}
 \le
 \frac{(S/n) \cdot c}{\vert \hat{\bm{M}}_2 \vert^{1/p}}. 
\end{align*}
Then we have $Vol(C_1^*)=Vol(C_2^*)=Vol(C_0)=Vol(C_3)$.
Denote $c^* = (S/n) \cdot c/\vert \hat{\bm{M}}_1 \vert^{1/p}$.
Figure \ref{figure:prob^*-1} represents the behavior of the coverage
probabilities of $C_0$, $C_1^*$, $C_2^*$ and $C_3$. 
The behavior of the coverage probabilities of $C_1^*$ and $C_2^*$ are
almost the same and two dashed lines which represent the behavior of
them overlap each other and look like just one line.  
We can see from Figure \ref{figure:prob^*-1} that the coverage
probability of $C_1^*$ and $C_2^*$ are larger than that of $C_0$ and
larger or at least the same level as that of $C_3$.  
Especially when $\lambda$ is small, the coverage probability of $C_1^*$
and $C_2^*$ are considerably larger than that of $C_0$ and $C_3$.
Since the volumes of $C_1^*$ and $C_2^*$ are equal to 
$C_0$ and $C_3$, we can say that $C_1^*$ and $C_2^*$ are useful for 
practical use.

\section{Concluding remarks}
\label{sec:6}
In this article we considered the estimation of the MSE and the MSE matrix
of shrinkage estimators from a decision theoretical viewpoint.
We generalize the argument of Carter et al.\cite{Carter-etal}, 
Wan et al.\cite{Wan-etal} and Kubokawa and
Srivastava\cite{Kubokawa-Srivastava} and provided some nonnegative and
positive estimators improving on the UMVUE. 
We also proposed some confidence sets by using the positive
definite estimator of the MSE matrices and showed through Monte
Carlo studies that proposed confidence sets have equal volume and attain
higher coverage probabilities than the conventional one.  These results
suggest that the proposed estimators are useful for evaluating the
precision of shrinkage estimators. 

Carter et al.\cite{Carter-etal} and Wan et al.\cite{Wan-etal} derived 
asymptotic coverage probabilities of their proposed confidence region.
By using their argument, we may provide better confidence region.

The argument in this article may also be available to the estimation of
MSE and MSE matrix in the small-area problem discussed in Prasad and
Rao\cite{Prasad-Rao}. 
Since the small-area estimation is an illustrative example of the
application of shrinkage estimators(e.g. Fay and
Herriot\cite{Fay-Herriot}), it would be interesting if we could extend
the results here to the problem
We leave these problem to our future research.

\section*{Appendix}
\appendix
\section{The derivation of (\ref{eq:MSE_positive}) and
 (\ref{eq:matrix_positive})} 
We first derive (\ref{eq:MSE_positive}).
Following Efron and Morris(1976), 
\begin{eqnarray*}
  R(\bm{\delta}_{\bm{\phi}})
  &=&
  \sigma^2 
  \mathrm{E}
  \left[
  p-
  \left(
  2(p-2)\frac{\phi(W)}{W}-(n+2)\frac{\phi^2(W)}{W}
  +4\phi'(W)+4\phi(W)\phi'(W)
  \right)
  \right]\\
 &=&
 \sigma^2 \mathrm{E}\left[p-
  \left(
  2(p-2)\frac{\phi(W)}{W}+4\phi'(W)
  \right)
  \right]
  +
  \mathrm{E}\left[
  S\frac{\phi^2(W)}{W}
  \right].
\end{eqnarray*}
The UMVU estimator of 
\(
\sigma^2 \mathrm{E}[p]
\)
is
\(pS/n\).
With respect to the term
\[
\sigma^2 \mathrm{E}\left[
  2(p-2)\frac{\phi(W)}{W}+4\phi'(W)
  \right], 
\]
we use the \(\chi^2\) identity with an absolutely continuous function 
\(g(W)\) by Efron and Morris(1976), 
$$
\mathrm{E}[g(W)S]
=
\sigma^2 \mathrm{E}[
ng(W) - 2g'(W)W].
$$
Let $g(W)$ be a solution of the following 
differential equation of first order,
$$
ng(W)-2g'(W)W = 2h(W),
$$
where
$$
h(W) \equiv (p-2)\frac{\phi(W)}{W} +2\phi'(W).
$$
Then Wan et al.\cite{Wan-etal} showed that the UMVUE of
$\hat{R}_0(\bdeltap)$ can be written as in (\ref{UMVUE:MSE}).
The general solution of $g(W)$ is written by 
\begin{equation}
 \label{solution:g(W)}
g(W)
=W^{\frac{n}{2}}
  \left\{
  \int_W^{\infty}
  t^{-\frac{n}{2}} \cdot
  \frac{h(t)}{t}dt
  + C_0 \right \}, 
\end{equation}
where \(C_0\) is a constant.
We note that $g(W)$ is required to be absolutely continuous.
When $h(W)$ is continuous, $g(W)$ is also continuous.
From the assumption 
$\phi(W) < \infty$ and $\phi'(W) < \infty$,
$\mathrm{E}[g(W)]$ exists if and only if $C_0=0$. 
However in the case of $\bdelta_{JS}^+$, 
$h(W)$ is not continuous and 
$g(W)$ with $C_0=0$ is also not continuous.
We can easily show that 
$g(W)$ is absolutely continuous and $E[g(W)] < \infty$ 
if and only if $C_0$ satisfies 
$$
C_0 = \left\{
\begin{array}{ll}
 \displaystyle{
2
\left(
\frac{p}{n} - \frac{p-2}{n+2}
\right)
\cdot
\left(
\frac{p-2}{n+2}
\right)^{-n/2}}, & \text{if } 
\displaystyle{W \le \frac{p-2}{n+2}}\\
0, & \text{otherwise}.
\end{array}
\right.
$$
By following Wan et al.\cite{Wan-etal} and the above argument,
the derivation of (\ref{eq:matrix_positive}) is similar and omitted. 
\hfill\qed

\section{Proof of Theorem \ref{th:MSE-2}}
The proof is similar to the proof of Theorem \ref{th:MSE}. 
Let $\psi^*(W)$ be the function which minimizes 
$\mathrm{E}[(\hat{R}^*(\psi;\bdeltap) - R^*(\bdeltap))^2 \mid W]$. 
Then 
\begin{align*}
 \psi^*_{\lambda}(W)
 &= \frac{n}{p} \cdot
 \frac{\mathrm{E}[S/\sigma^2 \mid W]}
 {\mathrm{E}[(S/\sigma^2)^2 \mid W]} 
 \cdot
 \frac{R^*(\bm{\delta}_{\bm{\phi}})}{\sigma^2}
 \cdot 
 \frac{1}{a(W)}\\
 & \le \frac{n}{p} \cdot \frac{1+W}{n+p+2} \cdot
 \frac{\alpha_{p,n}(\bdeltap)}{a(W)}\\
 & \le 
 \frac{n}{p} \cdot \frac{1+W_{p,n}(\bdeltap)}{n+p+2} \cdot
 \frac{\alpha_{p,n}(\bdeltap)}{a(W)}.
\end{align*}
By using the similar argument in the proof of Theorem \ref{th:MSE}, 
we can show that $\hat{R}^*(\psi_1;\bdeltap)$ dominates 
$\hat{R}^*_0(\bdeltap)$. 
\hfill\qed

\section{Proof of $\hat{R}^*(\psi_1;\bdelta_{JS}) >0$}
Since $\alpha_{p,n}(\bdelta_{JS}) < p$, 
it suffices to show that $W_{p,n}(\bdelta_{JS}) < (p+2)/n$.
In the case of $\bdelta_{JS}$, 
$a(W)=n(p-2)^2/(p(n+2)^2 W)$ from (\ref{eq:MSE_JS}) and 
$\alpha_{p,n}(\bdelta_{JS}) = (n(p-2))/(n+2)$(e.g. Kubokawa and
Srivastava\cite{Kubokawa-Srivastava}). 
Hence (\ref{eq:W}) is equivalent to 
$$
\zeta(W) := W(1+W) = 
\frac{n+p+2}{n}\cdot
\frac{p-2}{n+2}.
$$
Since $\zeta(W)$ is increasing for $W \ge 0$ and 
$$
\zeta\left(
\frac{p+2}{n} 
\right)=
\frac{p+2}{n} \cdot 
\frac{n+p+2}{n} \ge 
\frac{n+p+2}{n}\cdot
\frac{p-2}{n+2}, 
$$
we have $W_{p,n}(\bdelta_{JS}) < (p+2)/n$.
\hfill\qed

\section{Proof of Theorem \ref{th:MSE-3}}
For an estimator $\hat{R}^*(\delta_{\phi})$, 
denote
$
D(\hat{R}^*(\delta_{\phi}))
=
\mathrm{E}[
L(\hat{R}^*(\bdeltap) ; R^*(\bdeltap))
]
$.
It suffices to show that when $\psi(\cdot)$ satisfies the condition of
the theorem, 
$D(\hat{R}^*_0(\delta_{\phi})) - 
 D(\hat{R}^*(\psi;\delta_{\phi})) \ge 0$.
By using the same argument as the proof of Theorem 3 in Kubokawa and 
Srivastava\cite{Kubokawa-Srivastava}, we have
{\allowdisplaybreaks
\begin{align}
 \label{IERD:MSE}
 & D(\hat{R}^*_0(\delta_{\phi})) - 
 D(\hat{R}^*(\psi;\delta_{\phi}))\notag\\
 & \quad =
 \mathrm{E}
 \Biggl[
 \int_{1}^{\infty}
 \frac{d}{dt}
 \biggl\{
 \left(
 \frac{pS}{n}
 \right)^2
 a^2(W) \psi^2(tW)  -
 2  \left(
 \frac{pS}{n}
 \right)
 a(W) \psi(tW) R^*(\bdeltap)
 \biggr\} dt
 \Biggr]\notag\\
 & \quad \ge
 \mathrm{E}
 \Biggl[
 \int_{1}^{\infty}
 \bigg\{
 2 \left(
 \frac{pS}{n}
 \right)^2
 a^2(W) \psi'(tW) -
 \left(
 \frac{pS}{n}
 \right)
 a(W) \alpha_{p,n}(\bdeltap)
 \biggr\}W \psi'(W)dt
 \Biggr]\notag\\
 &\quad = 
 \frac{2p}{n}
 \times
 \int\int\int_0^{vw}
 \left\{
 \frac{p}{n} \cdot v a^2 (x/v)\psi(w) - 
 a(x/v) \alpha_{p,n}(\bdeltap)
 \right\}\notag\\
 & \qquad\qquad\qquad\qquad \times
 \psi'(w) f_p(x;\lambda) f_n(v) dxdvdw\notag\\
 &\quad =\int\mathrm{E}^v
 \Biggl[
 \int_0^{vw}  g_1(x/v)
 f_p(x;\lambda) dx \notag\\
 &\qquad\qquad\qquad\qquad
 \times
 \left\{
 \frac{p}{n}\cdot
 \frac{\int_0^{vw} a^2(x/v)f_p(x;\lambda)dx}
 {\int_0^{vw} a(x/v)f_p(x;\lambda)dx} \cdot
 v\psi(w) - \alpha_{p,n}(\bdeltap)
 \right\}
 \Biggr]\phi'(w)dw\notag\\
 &\quad \ge
 \frac{2p}{n}
 \int\mathrm{E}^v
 \Biggl[
 \int_0^{vw} a(x/v) f(x;\lambda) dz
 \cdot 
 \biggl\{
 \frac{p}{n}\cdot
 v a(w) \psi(w)
 - \alpha_{p,n}(\bdeltap)
 \biggr\}
 \Biggr]\phi'(w)dw\notag\\
 &\quad \ge
 \inf_w
 \Biggl\{
 \frac{2p}{n}\cdot
 \mathrm{E}^v
 \Biggl[
 \int_0^{vw} a(x/v) f_p(x;\lambda) dz
 \cdot 
 \biggl\{
 \frac{p}{n}\cdot
 v a(w) \psi(w)
 - \alpha_{p,n}(\bdeltap)
 \biggr\}
 \Biggr]\phi'(w)\Biggr\}, 
\end{align}
}
where $\mathrm{E}^v[\cdot]$ denotes the expectation with respect to 
$v \sim \chi^2_p$.
The first inequality follows from the fact that 
$\alpha_{p,n}(\bdeltap) \ge R^*(\bdeltap)$.
The second equality follows from the transformation, 
$(x,v,w)=(WS,S,tW)$ with $dt dW dS = x^{-1} dx dv dw$.
The second inequality follows from 
the nonincreasingness of $a(\cdot)$.
Since $a(\cdot)$ is nonincreasing, 
$\int_0^{vw} a(x/v) f(x;\lambda) dz$ 
is nondecreasing in $v$.
With respect to 
$q(v) = (p/n)va(w)\phi(w) -
\alpha_{p,n}(\bdeltap)$,  
there exists $v_0$ such that 
$q(v) \le 0$ for $0 < v \le v_0$ and 
$q(v) > 0$ for $v > v_0$.
Then we have from Kubokawa\cite{Kubokawa-1998}
\begin{align*}
 & \mathrm{E}^v
 \Biggl[
 \int_0^{vw} a(x/v) f_p(x;\lambda) dz
 \cdot 
 \biggl\{
 \frac{p}{n}\cdot
 v a(w) \psi(w)
 - \alpha_{p,n}(\bdeltap)
 \biggr\}
 \Biggr]\\
 &\qquad \ge 
 \int_0^{v_0w} a(x/v) f_p(x;\lambda) dz
 \cdot 
 \mathrm{E}^v
 \Bigl[
 \frac{p}{n}\cdot
 v a(w) \psi(w)
 - \alpha_{p,n}(\bdeltap)
 \Bigr].
\end{align*}
If $\psi(\cdot)$ satisfy the condition of the theorem, 
$$
\mathrm{E}^v
\left[
\frac{p}{n}\cdot
v a(w) \psi(w) - \alpha_{p,n}(\bdeltap)
\right]
=
p a(w) \psi(w)- \alpha_{p,n}(\bdeltap)
\ge 0.
$$ 
Therefore we can complete the proof.
\hfill\qed

\section{Proof of Theorem \ref{th:matrix-1}}
Define $\bm{M}_I(\bdeltap) = \sigma^{-2}\bm{M}(\bdeltap)$ and 
$\bm{M}^*_I(\bdeltap) = \bm{I}_p - \sigma^{-2} \bm{M}(\bdeltap)$.
Before we give a proof, we present a lemma required to prove the
theorem. 
\begin{lemma}
 \label{lemma:1}
 $\beta^{(1)}_{p,n}$ satisfies
 $$
 \beta^{(1)}_{p,n}(\bdeltap) \le 
 \mathrm{tr}
 \left(
 \bm{I}_p - 
 \frac{\bm{X}\bm{X}}{\Vert \bm{X} \Vert^2} 
 \right) \cdot
 \bm{M}^*_I(\bdeltap).
 $$
\end{lemma}

By following Section 2.2 in Kubokawa and
Srivastava\cite{Kubokawa-Srivastava},  
the proof of this lemma is easy and omitted here.

\ \\
{\it Proof of Theorem \ref{th:matrix-1}.}\; 
Let 
$\xi^*(\cdot)$ and $\eta^*(\cdot)$ be the functions 
$\xi(\cdot)$ and $\eta(\cdot)$ which minimize
$
\mathrm{E}[\mathrm{tr}
(\hat{\bm{M}}(\xi,\eta;\bm{\delta}_\phi) - M(\bm{\delta}_\phi))^2 \mid W]
$. 
Similar to (\ref{psi^*}), we have 
\begin{align*}
 \xi^*(W) &= \frac{1}{g_1(W)}
 \left(
 \frac{1}{n}
 -
 \frac{1}{\mathrm{E}[(S/\sigma^2)^2|W]}
 \cdot
 \mathrm{tr}\left\{
 \mathrm{E}\left[
 \frac{S}{\sigma^2}
 \left(
 \bm{I}_p - 
 \frac{\bm{X}\bm{X}'}{\Vert \bm{X} \Vert^2}
 \right) 
 \cdot \frac{\bm{M}_I(\bdeltap)}{p-1}
 \mid W
 \right]\right\}
 \right), \\
 & = \frac{1}{g_1(W)}
 \biggl(
 \frac{1}{n}
 -
 \frac{\mathrm{E}[S/\sigma^2|W]}{\mathrm{E}[(S/\sigma^2)^2|W]}\\
 & \qquad\qquad\qquad
 +
 \frac{1}{\mathrm{E}[(S/\sigma^2)^2|W]}
 \cdot
 \mathrm{E}\left[
 \frac{S}{\sigma^2}
 \mathrm{tr}
 \left(
 \bm{I}_p - 
 \frac{\bm{X}\bm{X}'}{\Vert \bm{X} \Vert^2}
 \right) 
 \cdot \frac{\bm{M}^*_I(\bdeltap)}{p-1}
 \mid W
 \right]
 \biggr), \\
\end{align*}
\begin{align*}
 \eta^*(W) = 
 \frac{1}{g_1(W)}
 \biggl(
 \frac{1}{n} 
 & + g_2(W) + 
 \frac{\phi^2(W)}{W} \notag\\
 &- 
 \left.
 \frac{1}{\mathrm{E}[(S/\sigma^2)^2|W]} \cdot
 \mathrm{tr}\left\{
 \mathrm{E}\left[
 \frac{S}{\sigma^2} \cdot 
 \frac{\bm{X}\bm{X}'}{\Vert \bm{X} \Vert^2}
 \cdot \bm{M}_I(\bdeltap)
 \mid W
 \right]\right\}
 \right).
\end{align*}
We note that $\xi^*(\cdot)$ does not depend on $\eta(\cdot)$ and that  
$\eta^*(\cdot)$ does not depend on $\xi(\cdot)$.
From (\ref{ineq:chi^2}), Lemma \ref{lemma:1},
the assumption that $\beta_{p,n}^{(1)}(\bdeltap) \ge 0$  
and the facts that  
\begin{equation}
 \label{eq:p-definite}
  \bm{I}_p - \frac{\bm{X}\bm{X}'}{W} 
  =\bm{\Gamma}'
  \left(
   \bm{I}_p - \bm{E}_{11}
  \right)\bm{\Gamma}
  \ge \bm{0}, \quad 
  \frac{\bm{X}\bm{X}'}{W} 
  \ge \bm{0}, \quad 
  M_I(\bdeltap) \ge \bm{0}.
\end{equation}
we have 
$$
\frac{1}{g_1(W)}
\left(
\frac{1}{n} -
\frac{1+W}{n+p+2}
\right) \le 
\psi^*(W)
\le 
\frac{1}{ng_1(W)}, 
$$
$$
\eta^*(W) \le 
\frac{1}{g_1(W)}
\left(
\frac{1}{n} + g_3(W)
\right). 
$$
By using the same argument as in the proof of Theorem \ref{th:MSE}, 
we can complete the proof.
\hfill\qed\\

\section{Proof of Theorem \ref{th:matrix-2}}
We present a lemma required to prove the theorem. 
\begin{lemma}
 \label{lemma:2}
 $\beta^{(2)}_{p,n}$ satisfies
 $$
 \beta^{(2)}_{p,n}(\bdeltap) \ge 
 \mathrm{tr}
 \frac{\bm{X}\bm{X}}{\Vert \bm{X} \Vert^2} \cdot
 \bm{M}^*_I(\bdeltap), 
 $$
 $$
 (p-1)\beta^{(2)}_{p,n}(\bdeltap) \ge 
 \mathrm{tr}
 \left(
 \bm{I}_p - 
 \frac{\bm{X}\bm{X}}{\Vert \bm{X} \Vert^2} 
 \right) \cdot
 \bm{M}^*_I(\bdeltap).
 $$
\end{lemma}

By following Kubokawa and Srivastava\cite{Kubokawa-Srivastava}, 
the proof of this lemma is also easy and omitted here.

\ \\
{\it Proof of Theorem \ref{th:matrix-2}.}\; 
The proof is similar to the one of Theorem \ref{th:matrix-1}.
Let 
$\xi^*(\cdot)$ and $\eta^*(\cdot)$ be the functions 
$\xi(\cdot)$ and $\eta(\cdot)$ which minimize
$\mathrm{E}[\mathrm{tr}
(\hat{\bm{M}}^*(\xi,\eta;\bm{\delta}_\phi) - 
\bm{M}^*(\bm{\delta}_\phi))^2 \mid W]$. 
Then we have 
\begin{align*}
 \xi^*(W) & = \frac{1}{g_1(W)} \cdot
 \frac{1}{\mathrm{E}[(S/\sigma^2)^2|W]}
 \cdot
 \mathrm{tr}\left\{
 \mathrm{E}\left[
 \frac{S}{\sigma^2}
 \left(
 \bm{I}_p - 
 \frac{\bm{X}\bm{X}'}{W}
 \right) 
 \cdot \frac{\bm{M}^*_I(\bdeltap)}{p-1}
 \mid W
 \right]\right\}, 
\end{align*}
\begin{align*}
 \eta^*(W) &= 
 \frac{1}{g_1(W)}
 \left(
 g_3(W) + 
 \frac{1}{\mathrm{E}[(S/\sigma^2)^2|W]} \cdot
 \mathrm{tr}\left\{
 \mathrm{E}\left[
 \frac{S}{\sigma^2} \cdot 
 \frac{\bm{X}\bm{X}'}{W}
 \cdot \bm{M}^*_I(\bdeltap)
 \mid W
 \right]\right\}
 \right).
\end{align*}
From (\ref{ineq:chi^2}), (\ref{eq:p-definite}) and Lemma \ref{lemma:2},
we have  
$$
\xi^*(W) \le 
\frac{1+W_{p,n}^{\xi}(\bdeltap)}{n+p+2} \cdot 
\frac{\beta^{(2)}_{p,n}}{g_1(W)}, \quad 
\eta^*(W) \le 
\frac{1}{g_1(W)}
\left(
g_3(W) + \frac{1+W_{p,n}^{\eta}(\bdeltap)}{n+p+2} 
\cdot \beta^{(2)}_{p,n}
\right).
$$
Hence by using the same argument as in the proof of Theorem
\ref{th:MSE}, we can complete the proof.
\hfill\qed

\section{Proof of Theorem \ref{th:matrix-3}}
For an estimator $\hat{\bm{M}}^*(\bdeltap)$, 
denote
$
D(\hat{\bm{M}}^*(\bdeltap)) =
\mathrm{E}[
L(\hat{\bm{M}}^*(\bdeltap) ; 
\bm{M}^*(\bdeltap))]
$.
Denote $\hat{\bm{M}}^*(\xi,\eta;\bdeltap)$ with $\eta(W)=1$
by $\hat{\bm{M}}^*(\xi;\bdeltap)$.
We first show that when $\xi(\cdot)$ satisfies the conditions of the 
theorem, $\hat{\bm{M}}^*(\xi;\bdeltap)$ dominates 
$\hat{\bm{M}}^*_0(\bdeltap)$. 
In the same way as (\ref{IERD:MSE}), we have 
{\allowdisplaybreaks
\begin{align}
 \label{IERD:matrix}
 & D(\hat{\bm{M}}^*_0(\bdeltap)) - 
 D(\hat{\bm{M}}^*(\xi;\bdeltap))\notag\\
 & \quad =
 \mathrm{E}
 \Biggl[
 \int_{1}^{\infty}
 \frac{d}{dt}
 \big\{
 (p-1)S^2 g_1^2(W) \xi^2(tW) \notag\\
 & \qquad\qquad\qquad -
 2 S g_1(W) \xi(tW) 
 \mathrm{tr}\left[
 \bm{\Gamma}'(\bm{I}_p - \bm{E}_{11})\bm{\Gamma}
 \bm{M}^*(\bdeltap)
 \right]
 \bigr\} dt
 \Biggr]\notag\\
 & \quad \ge
 (p-1)\sigma^4
 \mathrm{E}
 \Biggl[
 \int_{1}^{\infty}
 \frac{d}{dt}
 \big\{
 (S/\sigma^2)^2 g_1^2(W) \xi^2(tW) -
 2 (S/\sigma^2) g_1(W) \xi(tW) 
 \beta^{(2)}_{p,n}(\bdeltap)
 \bigr\} dt
 \Biggr]\notag\\
 & \quad =
 2(p-1)\sigma^4
 \mathrm{E}
 \Biggl[
 \int_{1}^{\infty}
 \big\{
 (S/\sigma^2)^2 g_1^2(W) \xi(tW) -
 (S/\sigma^2) g_1(W) \beta^{(2)}_{p,n}(\bdeltap)
 \bigr\}
 \xi'(W)dt
 \Biggr]\notag\\
 &\quad =\int\mathrm{E}^v
 \Biggl[
 \int_0^{vw} g_1(x/v)
 f_p(x;\lambda) dx \notag\\
 &\qquad\qquad\qquad\qquad
 \times
 \left\{
 \frac{\int_0^w g_1^2(x/v)f_p(x;\lambda)dx}
 {\int_0^w g_1(x/v)f_p(x;\lambda)dx} 
 \cdot v^2 \xi(w)
 - v \beta^{(2)}_{p,n}(\bdeltap)
 \right\}
 \Biggr]\xi'(w)dw\notag\\
 &\quad \ge
 \int\mathrm{E}^v
 \Biggl[
 \int_0^w g_1(x/v) f_p(x;\lambda) dz
 \cdot 
 \bigl\{
 v^2 g_1(w)\xi(w)
 - v \beta^{(2)}_{p,n}(\bdeltap)
 \bigr\}
 \Biggr]\xi'(w)dw.
\end{align}
}
The first inequality follows from Lemma \ref{lemma:2}
and the second inequality follows from the nonincreasingness of 
$g(\cdot)$. 
Similar to the argument in the proof of Theorem \ref{th:MSE-2},  
if $\xi(W)$ satisfy the conditions of the theorem, we have 
$$
\mathrm{E}^v[v^2 g_1(w)\xi(w) - v \beta^{(2)}_{p,n}(\bdeltap)] 
=
n(n+2)g_1(w)\xi(w) - n \beta^{(2)}_{p,n}(\bdeltap)
\ge 0, 
$$ 
and hence $\hat{\bm{M}}^*(\xi;\bdeltap)$ dominates 
$\hat{\bm{M}}^*_0(\bdeltap)$. 

Next we show that when $\eta(\cdot)$ satisfies the conditions of the
theorem, $\hat{\bm{M}}^*(\xi,\eta;\bdeltap)$ dominates 
$\hat{\bm{M}}^*(\xi;\bdeltap)$. 
In the same way as the above argument, we have 
{\allowdisplaybreaks
\begin{align}
 \label{IERD:matrix-2}
 & D(\hat{\bm{M}}^*(\xi;\bdeltap)) - 
 D(\hat{\bm{M}}^*(\xi,\eta;\bdeltap))\notag\\
 & \quad =
 \mathrm{E}
 \Biggl[
 \int_{1}^{\infty}
 \frac{d}{dt}
 \big\{
 S^2 g_1^2(W) \eta^2(tW) 
 - 2 S^2 g_1(W)
 \left(
 g_3(W) + \frac{\phi^2(W)}{W}
 \right) \eta(tW) 
 \notag\\
 & \qquad\qquad\qquad -
 2 S g_1(W) 
 \eta(tW) 
 \mathrm{tr}[
 \bm{\Gamma}'\bm{E}_{11}\bm{\Gamma}
 \bm{M}^*(\bdeltap)]
 \bigr\} dt
 \Biggr]\notag\\
 & \quad \ge
 \sigma^4
 \mathrm{E}
 \Biggl[
 \int_{1}^{\infty}
 \bigg\{
 (S/\sigma^2)^2 g_1^2(W) \eta^2(tW) -
 2 (S/\sigma^2)^2 g_1(W) 
 \left(
 g_3(W) + \frac{\phi^2(W)}{W}
 \right) \eta(tW) 
 \notag\\
 & \qquad\qquad\qquad -
 2 (S/\sigma^2)^2 g_1(W) \eta(tW) \beta^{(2)}_{p,n}(\bdeltap)
 \biggr\} dt
 \Biggr]\notag\\
 & \quad =
 2 \sigma^4
 \mathrm{E}
 \Biggl[
 \int_{1}^{\infty}
 \biggl\{
 (S/\sigma^2)^2 g_1^2(W) \eta(tW) -
 (S/\sigma^2)^2 g_1(W) 
 \left(
 g_3(W) + \frac{\phi^2(W)}{W}
 \right) 
 \notag\\
 & \qquad\qquad\qquad\qquad -
 (S/\sigma^2) g_1(W) \beta^{(2)}_{p,n}(\bdeltap)
 \biggr\} W \eta'(tW) dt
 \Biggr]\notag\\
 & \quad = 
 2 \sigma^4
 \int_0^{\infty}
 \mathrm{E}^v
 \Biggl[
 \int_{0}^{vw}
 \biggl\{
 v^2 g_1^2(z/v) \eta(w) \notag\\
 &\qquad\qquad\qquad\qquad -
 g^2_1(z/v) \cdot 
 \frac{v^2 g_3(z/v) + v \beta^{(2)}_{p,n}(\bdeltap)}{g_1(z/v)}
 \biggr\} f_p(z;\lambda) dz \Biggr]\eta'(w) dw\notag\\
 & \quad \ge
 2 \sigma^4
 \int_0^{\infty}
 \mathrm{E}^v
 \Biggl[
 \int_{0}^{vw}
 g_1^2(z/v)
 f_p(z;\lambda)dz\notag\\
 &\qquad\qquad\qquad\qquad\qquad
 \times
 \biggl\{
 v^2 \eta(w) -
 \frac{v^2 g_3(w) + v \beta^{(2)}_{p,n}(\bdeltap)}{g_1(w)}
 \biggr\} 
 \Biggr]\eta'(w) dw
\end{align}
}
The first inequality follows from Lemma \ref{lemma:2}
and the second inequality follows from the nonincreasingness of 
$g_1(W)$ and the nondecreasingness of $g_3(W)/g_1(W)$. 
By using the same argument in the proof of Theorem \ref{th:MSE-3}, 
if $\eta(w)$ satisfies 
$$
\eta(w) \ge 
\min \left(
1,\; 
\frac{g_3(w)}{g_1(w)} 
+ 
\frac{1}{n+2} \cdot \frac{\beta^{(2)}_{p,n}(\bdeltap)}{g_1(w)}
\right), 
$$
the right hand side of (\ref{IERD:matrix-2}) is nonnegative.

\hfill\qed

\bibliographystyle{hplain}
\bibliography{MSE}

\end{document}